\documentclass{amsart}

\usepackage{amssymb, amsmath}
\usepackage{mathrsfs}
\usepackage{amscd}
\usepackage{verbatim}

\usepackage[colorlinks,linkcolor={blue},
citecolor={blue},urlcolor={red},]{hyperref}

\theoremstyle{plain}
\newtheorem{theorem}{Theorem}[section]
\theoremstyle{remark}
\newtheorem{remark}[theorem]{Remark}
\newtheorem{example}[theorem]{Example}
\theoremstyle{plain}
\newtheorem{corollary}[theorem]{Corollary}
\newtheorem{lemma}[theorem]{Lemma}
\newtheorem{proposition}[theorem]{Proposition}

\numberwithin{equation}{section}

\def\Z{{\mathbb Z}}

\def\R{{\mathbb R}}
\def\C{{\mathbb C}}

\newcommand{\E}{{\mathbb E}}
\renewcommand{\P}{{\mathbb P}}
\newcommand{\A}{{\mathcal A}}
\newcommand{\F}{{\mathcal F}}

\newcommand{\g}{\gamma}

\newcommand{\e}{\varepsilon}

\renewcommand{\O}{\Omega}

\newcommand{\calL}{{\mathcal B}}
\newcommand{\n}{\Vert}
\newcommand{\one}{{{\bf 1}}}

\newcommand{\lb}{\langle}
\newcommand{\rb}{\rangle}

\begin{document}

\title
[$R$-boundedness of smooth operator-valued
functions]{$R$-boundedness of smooth operator-valued functions}

\author{Tuomas Hyt\"onen}
\address{Department of Mathematics and Statistics\\
University of Helsinki\\ Gustaf H\"all\-str\"omin katu 2B \\ FI-00014 Helsinki\\ Finland}
\email{tuomas.hytonen@helsinki.fi}
\thanks{Tuomas Hyt\"onen is supported by the Academy of Finland (grant 114374).}

\author{Mark Veraar}
\address{Delft Institute of Applied Mathematics\\
Delft University of Technology \\ P.O. Box 5031\\ 2600 GA Delft\\The
Netherlands} \email{m.c.veraar@tudelft.nl, mark@profsonline.nl}
\thanks{Mark Veraar is supported by the Alexander von Humboldt foundation.
His visit to Helsinki, which started this project,
was funded by the Finnish Centre of Excellence in Analysis and Dynamics Research.}




\begin{abstract}
In this paper we study $R$-boundedness of operator families $\mathcal{T}\subset
\calL(X,Y)$, where $X$ and $Y$ are Banach spaces. Under cotype and type
assumptions on $X$ and $Y$ we give sufficient conditions for $R$-boundedness.
In the first part we show that certain integral operator are $R$-bounded. This
will be used to obtain $R$-boundedness in the case that $\mathcal{T}$ is the
range of an operator-valued function $T:\R^d\to \calL(X,Y)$ which is in a
certain Besov space $B^{d/r}_{r,1}(\R^d;\calL(X,Y))$. The results will be
applied to obtain $R$-boundedness of semigroups and evolution families, and to
obtain sufficient conditions for existence of solutions for stochastic Cauchy
problems.
\end{abstract}

\maketitle

\section{Introduction}\label{sec:intro}

The notion of $R$-boundedness (see Section~\ref{ss:Rbd} for definition) appeared implicitly in the work of
Bourgain~\cite{Bou2} and was formalized by Berkson and Gillespie~\cite{BG}.
Cl\'ement, de Pagter, Suckochev and Witvliet~\cite{CPSW} studied it in more
detail in relation to vector-valued Schauder decompositions, and shortly after
Weis~\cite{We} found a characterization of maximal regularity for the Cauchy
problem $u'=Au+f$, $u(0)=0$, in terms of $R$-boundedness of the resolvent of
$A$ or the associated semigroup. After this, many authors have used
$R$-boundedness techniques in the theory of Fourier multipliers and Cauchy
problems (cf.\ \cite{DHP,Hyt04,KuWe} and references therein).

For Hilbert space operators, $R$-boundedness is equivalent to uniform boundedness. The basic
philosophy underlying much of the work cited above is that many results for Hilbert spaces remain true in certain
Banach spaces if one replaces boundedness by $R$-boundedness. Thus it is useful
to be able to recognize $R$-bounded sets of operators.

Let $X$ and $Y$ be Banach spaces. In this paper we will study
$R$-boundedness of some subsets of $\calL(X,Y)$ under type and
cotype assumptions. Although the definition of $R$-boundedness
suggests connections with type and cotype, there are
only few results on this in the literature. Arendt and Bu \cite[Proposition 1.13]{ArBu}
pointed out that uniform boundedness already implies
$R$-boundedness if (and only if) $X$ has cotype~$2$
and $Y$ has type~$2$. Recently, van Gaans~\cite{vG} showed that a countable union of $R$-bounded sets remains $R$-bounded if the individual $R$-bounds are $\ell^r$ summable for an appropriate $r$ depending on the type and cotype assumptions, improving on the trivial result with $r=1$ (the triangle inequality!) valid for any Banach spaces. Implicitly, one can find similar ideas already in Figiel~\cite{Figiel}.


In \cite[Theorem 5.1]{GW}, Girardi and Weis have found criteria for
$R$-boundedness of the range of operator-valued functions $T:\R^d\to
\calL(X,Y)$ in terms of their smoothness and the Fourier type of the Banach
space $Y$. Their result states that if $Y$ has Fourier type $p\in [1, 2]$ and
$T$ is in the Besov space $B^{\frac{d}{p}}_{p,1}(\R^d;\calL(X,Y))$, then
$\{T(t):t\in \R^d\}$ is $R$-bounded.

We will prove a similar result as \cite[Theorem 5.1]{GW} under
assumptions on the cotype of $X$ and the type of $Y$. More
precisely, if $X$ has cotype $q$ and $Y$ has type $p$ and if $T\in
B^{\frac{d}{r}}_{r,1}(\R^d;\calL(X,Y))$ for some $r\in [1, \infty]$
such that $\frac1r =\frac1p-\frac1q$, then $\{T(t):t\in \R^d\}$ is
$R$-bounded (see Theorem \ref{thm:main} below). Our result improves
\cite[Theorem 5.1]{GW}. This follows from the fact that every space
with Fourier type~$p$ has type~$p$. Furthermore, we note that the
only spaces which have Fourier type~$2$ are spaces which are
isomorphic to a Hilbert space. However, there are many Banach spaces
with type~$2$, e.g., all $L^p$ spaces with $p\in [2, \infty)$. In
the limit case that $X$ has cotype $2$ and $Y$ has type $2$ our
assumption on $T$ becomes $T\in B^{0}_{\infty,1}(\R^d;\calL(X,Y))$.
This condition is quite close to uniform boundedness of $\{T(t):
t\in \R^d\}$ which under these assumption on $X$ and $Y$ is
equivalent to $R$-boundedness.

Following \cite[Section 5]{GW}, we apply the sufficient condition for
$R$-boundedness to strongly continuous semigroups. Furthermore, we show
that our results are sharp in the case of the translation semigroup on $L^p(\R)$. The
$R$-boundedness result for semigroups leads to existence, uniqueness and
regularity results for stochastic equations with additive noise. As a second
application we present an $R$-boundedness result for evolution families,
assuming the conditions of Acquistapace and Terreni~\cite{AT2}.

We will write $a\lesssim b$ if there exists a universal constant
$C>0$ such that $a\leq Cb$, and $a\eqsim b$ if $a\lesssim b\lesssim
a$. If the constant $C$ is allowed to depend on some parameter $t$,
we write $a\lesssim_t b$ and $a\eqsim_t b$ instead.

\section{Preliminaries}

Throughout this paper $(\O,\A,\P)$ denotes a probability space, and $\E$ is the expectation. Let
$X$ and $Y$ be Banach spaces. Let $(r_n)_{n\geq 1}$ be a Rademacher
sequence on $\O$, i.e.\ an independent sequence with
\[\P(r_n = 1) = \P(r_n = -1) = \frac12.\]
For $N\geq 1$ and $x_1, \ldots, x_N\in E$, recall the
Kahane-Khinchine inequalities (cf.\ \cite[Section 11.1]{DJT} or
\cite[Proposition 3.4.1]{KwWo}): for all $p,q\in [1, \infty)$, we
have
\begin{equation}\label{eq:KKineq}
\Big(\E\Big\|\sum_{n=1}^N r_n x_n\Big\|^p\Big)^{\frac1p} \eqsim_{p,q}
\Big(\E\Big\|\sum_{n=1}^N r_n x_n\Big\|^q\Big)^{\frac1q}.
\end{equation}
These inequalities will be often applied without referring to it explicitly.

For each integer $N$, the space ${\rm Rad}_N(X)\subset L^2(\O;X)$ is defined as
all elements of the form $\sum_{n=1}^N r_n x_n$, where $(x_n)_{n=1}^N$ are in
$X$.

\subsection{Type and cotype\label{subsec:type}}

Let $p\in [1,2]$ and $q\in [2,\infty]$. The space $X$ is said to have {\em type
$p$} if there exists a constant $C\ge 0$ such that for all $(x_n)_{n=1}^N$ in
$X$ we have
$$
\Big(\E \Big\n \sum_{n=1}^N r_n x_n\Big\n^2\Big)^\frac12 \le C
\Big(\sum_{n=1}^N \n x_n\n^p\Big)^\frac1p.
$$
The space $X$ is said to have {\em cotype $q$} if there exists a constant $C\ge
0$ such that for all $(x_n)_{n=1}^N$ in $X$ we have
$$
\Big(\sum_{n=1}^N \n x_n\n^q\Big)^\frac1q \le C\Big(\E \Big\n \sum_{n=1}^N r_n
x_n\Big\n^2\Big)^\frac12,
$$
with the obvious modification in the case $q=\infty$.

For a detailed study of type and cotype we refer to \cite{DJT}. Every Banach
space has type $1$ and cotype $\infty$ with constant $1$. Therefore, we say
that $X$ has {\em non-trivial type (non-trivial cotype)} if $X$ has type $p$
for some $p\in (1, 2]$ (cotype $q$ for some $2\leq q<\infty$). If the space $X$
has non-trivial type, it has non-trivial cotype. Hilbert spaces have type $2$
and cotype $2$ with constants $1$. For $p\in [1, \infty)$ the $L^p$-spaces have
type $p\wedge 2$ and cotype $p\vee 2$.

Recall the following duality result for ${\rm Rad}_N(X)$ (cf. \cite{PisConv} or
\cite[Chapter 13]{DJT}). If $X$ has non-trivial type then
\begin{equation}\label{eq:Rad}
{\rm Rad}_N(X)^* = {\rm Rad}_N(X^*)
\end{equation}
isomorphically with constants independent of $N$.

\subsection{Fourier type}

The Fourier transform $\widehat{f}=\F f$ of a function $f\in L^1(\R^d;X)$ will be
normalized as
$$ \widehat{f}(\xi) = \frac1{(2\pi)^{d/2}}\int_{\R^d} f(x)e^{-ix\cdot\xi}\,dx, \quad \xi\in\R^d.$$

Let $p\in [1, 2]$ and $p'$ be the conjugate exponent, $\frac1p+\frac{1}{p'}=1$.
The space $X$ has {\em Fourier type $p$}, if $\F$ defines a bounded linear operator
 for some (and then for all) $d=1, 2, \ldots$ from $L^p(\R^d;X)$ to $L^{p'}(\R^d;X)$.

If $X$ has Fourier type $p$, then it has both type $p$ and cotype $p'$.
In particular, spaces isomorphic to a Hilbert space are the only ones with Fourier type~$2$ (see~\cite{Kwa}).
The $L^p$-spaces have Fourier type $p\wedge p'$ (see \cite{Pe69}), while
every Banach space has Fourier type $1$. The notion becomes more restrictive with increasing $p$.



\subsection{$R$-boundedness}\label{ss:Rbd}

A collection $\mathcal{T}\subset \calL(X,Y)$ is said to be {\em
$R$-bounded}\index{R-boundedness}\index{$R$-boundedness}\index{operator!$R$-bounded}
if there exists a constant $M\ge 0$ such that
\[
\Bigl(\E\Bigl\|\sum_{n=1}^N r_n T_n x_n \Bigr\|_{Y}^2\Bigr)^\frac12
\leq M\Bigl(\E\Bigl\|\sum_{n=1}^N r_n x_n
\Bigr\|_{X}^2\Bigr)^\frac12,
\]
for all $N\ge 1$ and all sequences $(T_n)_{n=1}^N$ in $\mathcal{T}$
and $(x_n)_{n=1}^N$ in $X$. The least constant $M$ for which this
estimate
holds is called the {\em $R$-bound}\label{p:rbound} of
$\mathcal{T}$, notation $R({\mathcal T})$. By \eqref{eq:KKineq}, the
role of the exponent $2$ may be replaced by any exponent $1\le p <
\infty$ (at the expense of a possibly different constant).

The notion of $R$-boundedness has played an important role in recent
progress in Fourier multiplier theory and this has applications to
regularity theory of parabolic evolution equations. For details on
the subject we refer to \cite{DHP,KuWe} and references therein.

A property which we will need later on is the following. If $\mathcal{T}\subset
\calL(X,Y)$ is $R$-bounded and $X$ has non-trivial type, then it follows from
\eqref{eq:Rad} that the set of adjoint operators $\mathcal{T}^* = \{T^*\in
\calL(Y^*,X^*): T\in \mathcal{T}\}$ is $R$-bounded as well.

\subsection{Lorentz spaces}
We recall the definition of the Lorentz space (cf. \cite{Gra,Tr1}).
Let $(S,\Sigma, \mu)$ be a $\sigma$-finite measure space. For $f\in
L^1(S)+L^\infty(S)$
define the non-increasing rearrangement of $f$ as
\[f^*(s) = \inf\{t>0: \mu(|f|>t)\leq s\}, \ s>0.\]
For $p,q\in [1,\infty]$ define
\[L^{p,q}(S) = \{f\in L^1(S)+L^\infty(S): \|f\|_{L^{p,q}(S)}<\infty\},\]
where
\begin{equation*}
  \|f\|_{L^{p,q}(S)}
  =\begin{cases}
    \Big(\int_0^\infty t^{q/p} f^*(t)^q\frac{dt}{t}\Big)^{1/q} & \text{if }q\in[1,\infty), \\
    \sup_{t>0}t^{1/p}f^*(t) & \text{if }q=\infty.
  \end{cases}
\end{equation*}
For $p\in [1,\infty]$ and $q_1<q_2$ one has
\[\|f\|_{L^{p,p}(S)} = \|f\|_{L^p(S)}, \ \  \|f\|_{L^{p,q_2}(S)}\leq c_{p,q_1,q_2}\|f\|_{L^{p,q_1}(S)}.\]

Also recall (e.g.~\cite[pp.~331--2]{Sawyer}) that for $p,q\in [1,\infty)$
\begin{equation}\label{eq:Lorentzdistr}
||f\|_{L^{p,q}(S)} = \Big(\int_0^\infty p\, t^q\big(\mu(|f|>t)\big)^{q/p}  \, \frac{dt}{t}\Big)^{\frac1q};
\end{equation}
indeed, just compare the two iterated integrals of $s^{q-1}t^{q/p-1}$ over the subset $\{f^*(t)>s\}=\{\mu(|f|>s)>t\}$ of $(0,\infty)^2$.

\subsection{Besov spaces\label{subsec:Besov}}
We recall the definition of Besov spaces using the so-called
Little\-wood-Paley decomposition (cf. \cite{BeLo,Tr1}). Let
$\phi\in{\mathscr S}(\R^d)$ be a fixed Schwartz function
whose Fourier transform $\widehat\phi$ is nonnegative and has
support in $\{\xi\in\R^d: \ \tfrac12\le |\xi|\le 2\}$ and which
satisfies
$$ \sum_{k\in\Z} \widehat\phi(2^{-k}\xi) =1 \quad\hbox{for $\xi\in \R^d\setminus\{0\}$}.$$
Such a function can easily be constructed (cf.\ \cite[Lemma
6.1.7]{BeLo}). Define the sequence $({\varphi_k})_{k\ge 0}$ in
${\mathscr S}(\R^d)$ by
$$\widehat{\varphi_k}(\xi) = \widehat\phi(2^{-k}\xi) \quad \text{for}\ \  k=1,2,\dots \quad \text{and} \ \ \widehat{\varphi_0}(\xi) = 1- \sum_{k\ge 1} \widehat{\varphi_k}(\xi), \quad \xi\in\R^d.$$

Similar as in the real case one can define ${\mathscr S}(\R^d;X)$ as
the usual {\em Schwartz space} of rapidly decreasing $X$-valued
smooth functions on $\R^d$. As in the real case this is a Fr\'echet
space. Let the space of $X$-valued {\em tempered
distributions}\index{tempered distributions} ${\mathscr S}'(\R^d;X)$
be defined as the continuous linear operators from ${\mathscr
S}(\R^d)$ into $X$.

For $1\le p,q \le \infty$ and $s\in\R$ the {\em Besov
space}\index{space!Besov}\index{Besov space} $B_{p,q}^s(\R^d;X)$ is
defined as the space of all $X$-valued tempered distributions $f\in
{\mathscr S}'(\R^d;X)$ for which
\[ \| f\|_{B_{p,q}^s (\R^d;X)} := \Big\| \big( 2^{ks}{\varphi}_k * f\big)_{k\ge 0}
\Big\|_{l^q(L^p(\R^d;X))}
\]
is finite. Endowed with this norm, $B_{p,q}^s(\R^d;X)$ is a Banach
space, and up to an equivalent norm this space is independent of the
choice of the initial function $\phi$. The sequence $({\varphi}_k *
f)_{k\ge 0}$ is called the {\em Littlewood-Paley decomposition} of
$f$ associated with the function $\phi$.

If $1\le p,q<\infty$, then $B_{p,q}^s(\R^d;X)$ contains the Schwartz
space ${\mathscr S}(\R^d;X)$ as a dense subspace.

For $1\leq p_1\leq p_2 \leq \infty$, $q\in [1, \infty]$ and $s_1,s_2\in \R$
with $s_1-\frac{d}{p_1} = s_2-\frac{d}{p_2}$ the following continuous inclusion
holds (cf. \cite[Theorem 2.8.1(a)]{Tr1}
\[B^{s_1}_{p_1,q}(\R^d;X)\hookrightarrow B^{s_2}_{p_2,q}(\R^d;X).\]

Next we give an alternative definition of Besov spaces. Let $I = (a,b)$ with
$-\infty\leq a< b\leq \infty$. For $h\in\R$ and a function $f:I\to X$ we define
the function $T(h)f:\R\to X$ as the translate of $f$ by $h$, i.e.
\[(T (h) \phi)(t) :=
  \begin{cases}
    \phi(t+h) & \text{if $t+h\in I$}, \\
    0 & \text{otherwise}.
  \end{cases}
\]
For $h\in \R$ put
\[ I[h] : =  \Big\{r\in I:\ r+  h \in I\Big\}.\]
For a strongly measurable function $f\in L^p(I;X)$ and $t>0$ let
\[\varrho_p(f,t) := \sup_{|h|\le t} \Big(\int_{I} \|T(h) f(r) - f(r)\|^p \, dr\Big)^{\frac1p}.\]
We use the obvious modification if $p=\infty$. For $p,q\in [1,\infty]$ and
$s\in (0,1)$ define
\[\Lambda^s_{p, q}(I;X) := \{f\in L^p(I;X): \|f\|_{\Lambda^s_{p, q}(I;X)}<\infty\},\]
where
\begin{equation}\label{besov}
\|f\|_{\Lambda^s_{p, q}(I;X)} = \Big(\int_I \|f(r)\|^p \, dr\Big)^{\frac1p} +
\Big(\int_0^1 \big(t^{-s}\varrho_p(f,t)\big)^q\, \frac{dt}{t}\Big)^\frac1q
\end{equation}
with the obvious modification if $p=\infty$ or $q=\infty$. Endowed with the
norm $\|\cdot\|_{\Lambda^s_{p, q}(I;X)}$, $\Lambda^s_{p,q}(I;X)$ is a Banach
space. Moreover, if $I=\R$, then $\Lambda^s_{p,q}(\R;X) = B^s_{p,q}(\R;X)$ with
equivalent norms (cf. \cite[Proposition 3.1]{PeWo} and \cite[Theorem
4.3.3]{Schm}). Similarly, if $I\neq \R$, then for every $f\in
\Lambda^s_{p,q}(I;X)$ there exists a function $g\in \Lambda^s_{p,q}(\R;X)$ such
that $g|_I = f$ and there exists a constant $C>0$ independent of $f$ and $g$
such that
\begin{equation}\label{eq:gf}
C^{-1}\|g\|_{\Lambda^s_{p,q}(\R;X)} \leq \|f\|_{\Lambda^s_{p,q}(I;X)} \leq
\|g\|_{\Lambda^s_{p,q}(\R;X)}.
\end{equation}

%

\section{Tensor products}

We start with a basic lemma, which can be viewed as a generalization
of the Kahane-contraction principle.

\begin{lemma}\label{lem:LqrboundedLorentz}
Let $X$ be a Banach space and let $(S,\Sigma, \mu)$ be a $\sigma$-finite
measure space and let $q\in [2,\infty)$. The following assertions hold:
\begin{enumerate}
\item\label{it:LorentzContraction} If $X$ has cotype $q$, then there exists a constant $C$ such that for all
$(f_n)_{n=1}^N$ in $L^{q,1}(S)$ and $(x_n)_{n=1}^N$ in $X$
\begin{equation}\label{eq:estLpnormfLorentz}
\begin{split}
  &\Big\|\sum_{n=1}^N r_n f_n x_n \Big\|_{L^q(S;L^2(\O;X))} \\
  &\leq C\int_0^{\infty}\max_{1\leq n\leq N}\mu(|f_n|>t)^{1/q}\,dt\,
    \Big\|\sum_{n=1}^N r_n x_n\Big\|_{L^2(\O;X)}.
\end{split}
\end{equation}
\item\label{it:LqContraction}
If $X$ has cotype $q$, then for all $\tilde{q}\in
(q,\infty]$ there exists a constant $C$ such that for all $(f_n)_{n=1}^N$ in
$L^{\tilde{q}}(S)$ and $(x_n)_{n=1}^N$ in $X$
\begin{equation}\label{eq:estLpnormf}
  \Big\|\sum_{n=1}^N r_n f_n x_n \Big\|_{L^{\tilde{q}}(S;L^2(\O;X))}
 \leq C \sup_{1\leq n\leq N}\|f_n\|_{L^{\tilde{q}}(S)} \Big\|\sum_{n=1}^N r_n x_n \Big\|_{L^2(\O;X)}.
\end{equation}
Moreover, if $q\in\{2,\infty\}$, then \eqref{eq:estLpnormf} holds with $\tilde{q}=q$.

\item\label{it:LorentzConverse} Assume $(S,\Sigma, \mu)$ contains $N$ disjoint sets of equal finite
positive measure  for every $N\in\Z_+$. If there exists a constant $C$ such
that \eqref{eq:estLpnormfLorentz} holds for all $(f_n)_{n=1}^N$ in $L^{q,1}(S)$
and $(x_n)_{n=1}^N$ in $X$, then $X$ has cotype~$q$.

\end{enumerate}
\end{lemma}

\begin{remark}
Note that by \eqref{eq:Lorentzdistr} for $q\in [1,\infty)$ it holds that
\[\int_0^{\infty}\max_{1\leq n\leq N}\mu(|f_n|>t)^{1/q}\,dt\leq \int_0^{\infty} \lambda\big(\max_{1\leq n\leq N} (f_n^*)>t\big)^{1/q}\,dt = \frac1q \big\|\max_{1\leq n\leq N} (f_n^*)\big\|_{L^{q,1}},\]
where $\lambda$ denotes the Lebesgue measure on $(0,\infty)$.
Moreover, if the $f_1, \ldots, f_N$ are identically distributed,
then one has
\[\int_0^{\infty}\max_{1\leq n\leq N}\mu(|f_n|>t)^{1/q}\,dt = \frac1q\|f_1\|_{L^{q,1}(S)}.\]

With this in mind, one can also view \eqref{eq:estLpnormfLorentz} as
an extension of \cite[Proposition 9.14]{LeTa} and \cite[Proposition
3.2(ii)]{Pi2}. There it is shown that \eqref{eq:estLpnormfLorentz}
holds for the case that $\mu(S) = 1$, $(f_n)_{n\geq 1}$ are i.i.d.\
and symmetric.
\end{remark}

\begin{remark}
By \eqref{eq:KKineq}, one could rephrase \eqref{eq:estLpnormf} as follows when
$\tilde{q}<\infty$. In the natural embedding $L^{\tilde{q}}(S)\hookrightarrow\calL(X,L^{\tilde{q}}(S;X))$,
$f\mapsto f\otimes(\cdot)$, the unit ball $B_{L^{\tilde{q}}(S)}$ becomes an $R$-bounded
subset of $\calL(X,L^{\tilde{q}}(S;X))$.
\end{remark}


\begin{proof}[Proof of Lemma \ref{lem:LqrboundedLorentz}]
\eqref{it:LorentzContraction} and \eqref{it:LqContraction}: \ 
Define the operator $T:\ell^\infty_N\to L^2(\O;X)$ by $T(a) = \sum_{n=1}^N r_n a_n x_n$.
By the Kahane contraction principle there holds
\[\|T\|_{\calL(\ell^\infty_N, L^2(\O;X))}\leq 2 \Big\|\sum_{n=1}^N r_n x_n
\Big\|_{L^2(\O;X)}.
\]
Since $L^2(\O;X)$ has cotype $q$ it follows from \cite[Theorem 11.14]{DJT} that
$T$ is $(q,1)$-summing with $\pi_{q,1}(T)\leq C_{X,q}\|T\|$. 
Then \cite[Theorem 2.1]{Pifact} (also see \cite[Theorem 10.9]{DJT})
implies that there is a probability measure $\nu$ on
$\{1,2,\ldots,N\}$ such that $T=\tilde{T}j$, where
$j:\ell^\infty_N\to \ell^{q,1}_N(\nu)$ is the embedding and
$\tilde{T}\in \calL(\ell^{q,1}_N(\nu),L^2(\O;X))$ satisfies
$\|\tilde{T}\|_{\calL(\ell^{q,1}_N(\nu),L^2(\O;X))}\leq
q^{-1+\frac1q}\pi_{q,1}(T)$. It follows that for all scalars
$(a_n)_{n=1}^N$,
\begin{equation*}
\begin{split}
  \Big\|\sum_{n=1}^N r_n a_n x_n \Big\|_{L^2(\O;X))}
  &\leq C_{X,q}\|(a_n)_{n=1}^N\|_{\ell^{q,1}_N(\nu)}
  \Big\|\sum_{n=1}^N r_n x_n \Big\|_{L^2(\O;X)} \\
  &\leq C_{X,q,\tilde{q}}\|(a_n)_{n=1}^N\|_{\ell^{\tilde{q}}_N(\nu)}
  \Big\|\sum_{n=1}^N r_n x_n \Big\|_{L^2(\O;X)},
\end{split}
\end{equation*}
where $\ell^{q,1}_N(\nu)$ denotes the Lorentz space $L^{q,1}$
defined on $\{1,\ldots, N\}$ with measure $\nu$, and the second step follows from the embedding of $\ell^{q,1}_N(\nu)$ into $\ell^{\tilde{q}}_N(\nu)$ for $\tilde{q}\in(q,\infty]$.

If we apply this with $a_n = f_n(s)$ and take the $L^q(\mu)$-norms, then it follows
from \eqref{eq:Lorentzdistr} and Minkowski's inequality
that
\begin{equation*}
\begin{split}
  &\Big\|\sum_{n=1}^N r_n f_n x_n \Big\|_{L^q(S;L^2(\O;X))}\times \Big\|\sum_{n=1}^N r_n x_n \Big\|_{L^2(\O;X)}^{-1}\\
  &\leq C_{X,q} \Big(\int_S \Big(\int_0^{\infty}\Big(\sum_{n=1}^N \nu(n)\one_{\{|f_n(s)|>t\}}
       \Big)^{1/q}\,dt \Big)^{q}\, d\mu(s) \Big)^{\frac1q}  \\
  &\leq C_{X,q} \int_0^{\infty}\Big(\int_S \sum_{n=1}^N \nu(n)\one_{\{|f_n(s)|>t\}} d\mu(s)\Big)^{1/q}\,dt \\
  &=C_{X,q}\int_0^{\infty}\Big(\sum_{n=1}^N\nu(n)\mu(|f_n|>t)\Big)^{1/q}\,dt \\
  &\leq C_{X,q}\int_0^{\infty}\Big(\max_{1\leq n\leq N}\mu(|f_n|>t)\Big)^{1/q}\,dt.
\end{split}
\end{equation*}
Similarly with $L^{\tilde{q}}(\mu)$-norms, it follows that
\begin{equation*}
\begin{split}
  &\Big\|\sum_{n=1}^N r_n f_n x_n \Big\|_{L^{\tilde{q}}(S;L^2(\O;X))}
     \times \Big\|\sum_{n=1}^N r_n x_n \Big\|_{L^2(\O;X)}^{-1}\\
  &\leq C_{X,q,\tilde{q}}\Big(\int_S \Big(\sum_{n=1}^N \nu(n)|f_n(s)|^{\tilde{q}}
       \Big)^{\tilde{q}/\tilde{q}}d\mu(s)\Big)^{1/\tilde{q}} \\
  &=C_{X,q,\tilde{q}}\Big(\sum_{n=1}^N\nu(n)\int_S|f_n(s)|^{\tilde{q}}\,d\mu(s)\Big)^{1/\tilde{q}} \\
  &\leq C_{X,q,\tilde{q}}\Big(\max_{1\leq n\leq N}\int_S |f_n(s)|^{\tilde{q}}\,d\mu(s)\Big)^{1/\tilde{q}}.
\end{split}
\end{equation*}

\eqref{it:LorentzConverse}: \ Let $x_1, \ldots, x_N\in X$. Let $(S_n)_{n=1}^N$ be disjoint
sets in $\Sigma$ with $\mu(S_n)=\mu(S_1)\in (0,\infty)$ for all $n$. Letting
$f_n = \mu(S_1)^{-1/q}\one_{S_n}$ for $n=1, 2, \ldots, N$, we
obtain that
\[\Big\|\sum_{n=1}^N r_n f_n x_n \Big\|_{L^q(S;L^2(\O;X))}^q = \sum_{n=1}^N \|x_n\|^q.\]
On the other hand,
\[
  \mu(|f_n|>t)=\mu(S_1)\cdot\one_{[0,\mu(S_1)^{-1/q})}(t)
\]
for all $n=1, \ldots,N$. Therefore, \eqref{eq:estLpnormfLorentz} implies that
\[\sum_{n=1}^N \|x_n\|^q \leq C^q \Big\|\sum_{n=1}^N r_n x_n
\Big\|_{L^2(\O;X)}^q,\] which shows that $X$ has cotype $q$.

\end{proof}

We do not know, whether we can take $\tilde{q}=q$ in Lemma~\ref{lem:LqrboundedLorentz}\eqref{it:LqContraction} if $q\neq 2$. However, if $X$ is an $L^{q}$-space with $q\geq 2$, then one may take $\tilde{q}=q$. This
follows from the next remark in the case that $\widetilde X =\R$.

\begin{remark}\label{rem:optimal}
Let $(A,\mathcal{A},\nu)$ be a $\sigma$-finite measure space. Let
$2\leq q_1<q<\infty$. Let $\widetilde{X}$ be a Banach space with
cotype $q_1$. If $X= L^q(A;\widetilde{X})$, then
\eqref{eq:estLpnormf} of
Lemma~\ref{lem:LqrboundedLorentz}\eqref{it:LqContraction} holds with
$\tilde{q}=q$.
\end{remark}
\begin{proof}
By Fubini's theorem, Lemma~\ref{lem:LqrboundedLorentz}\eqref{it:LqContraction} applied to
$\widetilde{X}$ and \eqref{eq:KKineq} we obtain that
\[\begin{aligned}
  \Big\|\sum_{n=1}^N r_n f_n x_n \Big\|_{L^q(S;L^2(\O;X))}
  & \eqsim_{q} \Big\|\sum_{n=1}^N r_n f_n x_n \Big\|_{L^{q}(A;L^{q}(S;L^2(\O;\widetilde{X})))} \\
  & \leq C \sup_{1\leq n\leq N}\|f_n\|_{L^{q}(S)} \Big\|\sum_{n=1}^N r_n x_n \Big\|_{L^{q}(A;L^2(\O;\widetilde{X}))}. \\
  & \eqsim_{q} C \sup_{1\leq n\leq N}\|f_n\|_{L^{q}(S)} \Big\|\sum_{n=1}^N r_n x_n \Big\|_{L^2(\O;X)}.
\end{aligned}
\]
\end{proof}

\begin{remark}
Notice that a version of Lemma \ref{lem:LqrboundedLorentz} also
holds for quasi-Banach spaces. This can be proved in a similar way
as above. Instead of \cite{Pifact} one has to use the factorization
result of \cite[Theorem 4.1]{KMS}. Note that in \cite{KMS} the role
of the Lorentz space $L^{q,1}(S)$ is replaced by $L^{q,r}(S)$, where
$r$ is some number in $(0,1]$ which depends on $X$. One can see from
the above proof that this number $r$ will also appear in the
quasi-Banach space version of \eqref{eq:estLpnormf}. The details are
left to the interested reader.
\end{remark}

The following dual version of Lemma~\ref{lem:LqrboundedLorentz} holds:

\begin{lemma}
Let $X$ be a Banach space, let $(S,\Sigma, \mu)$ be a
$\sigma$-finite measure space and let $p\in (1, 2]$. The following
assertions hold:
\begin{enumerate}
\item\label{it:dualLorContr} If $X$ has type $p$, then there exists a constant $C$ such that for all
$(f_n)_{n=1}^N$ in $L^{p,\infty}(S)$ which are identically
distributed and $(x_n)_{n=1}^N$ in $X$
\begin{equation}\label{eq:estLpnormfLorentzreverse}
\begin{split}
  &\|f_1\|_{L^{p,\infty}(S)} \Big\|\sum_{n=1}^N r_n x_n \Big\|_{L^2(\O;X)}
  \leq C\Big\|\sum_{n=1}^N r_n f_n x_n\Big\|_{L^p(S;L^2(\O;X))}.
\end{split}
\end{equation}

\item\label{it:dualLpContr} If $X$ has type $p$, then for all $\tilde{p}\in [1,p)$
there exists a constant $C$ such that for all $(f_n)_{n=1}^N$ in $L^{\tilde{p}}(S)$ and
$(x_n)_{n=1}^N$ in $X$
\begin{equation}\label{eq:estLpnormfreverse}
  \inf_{1\leq n\leq N} \|f_n\|_{L^{\tilde{p}}(S)} \Big\|\sum_{n=1}^N r_n x_n \Big\|_{L^2(\O;X)}
  \leq C \Big\|\sum_{n=1}^N r_n f_n x_n \Big\|_{L^{\tilde{p}}(S;L^2(\O;X))}.
\end{equation}
Moreover, if $p\in\{1,2\}$, then \eqref{eq:estLpnormf} holds with $\tilde{p}=p$.

\item\label{it:dualLorConverse} Assume $(S,\Sigma, \mu)$ contains $N$ disjoint sets of equal finite
positive measure  for every $N\in\Z_+$. If there exists a constant
$C$ such that \eqref{eq:estLpnormfLorentzreverse} holds for all
$(f_n)_{n=1}^N$ in $L^{p,1}(S)$ which are identically distributed
and $(x_n)_{n=1}^N$ in $X$, then $X$ has type~$p$.

\end{enumerate}
\end{lemma}

A similar statement as in Remark \ref{rem:optimal} also holds.

\begin{proof}
\eqref{it:dualLorContr} : Without loss of generality,
$\|f_1\|_{L^{p,\infty}(S)}=\sup_{t>0}t^{1/p}f_1^*(t)=1$. Choose
$t_0$ so that $f_1^*(t_0)>(2t_0)^{-1/p}$, or equivalently
$t_0<\mu\big(|f_1|>(2t_0)^{-1/p}\big)$. Let $A_n:=\{|f_n|>(2t_0)^{-1/p}\}$, so by
equidistribution, $\mu(A_n)=\mu(A_1)>t_0$. It follows that
\begin{equation*}
\begin{split}
   \Big| &\sum_{n=1}^N\langle x_n,x_n^*\rangle\Big|
   =\frac{1}{\mu(A_1)}\Big|\int_S\E\Big<\sum_{n=1}^N r_n \one_{A_n} x_n,
           \sum_{m=1}^N r_m \one_{A_m} x_m^*\Big>d\mu(s)\Big| \\
   &\leq\frac{1}{\mu(A_1)}\Big\|\sum_{n=1}^N r_n \one_{A_n} x_n\Big\|_{L^p(S;L^2(\O;X))}
        \Big\|\sum_{m=1}^N r_m \one_{A_m}
        x_m^*\Big\|_{L^{p'}(S;L^2(\O;X^*))}.
\end{split}
\end{equation*}

Now $X^*$ has cotype $p'$, hence by Lemma~\ref{lem:LqrboundedLorentz}\eqref{it:LorentzContraction} there
holds
\begin{equation*}
  \Big\|\sum_{m=1}^N r_m \one_{A_m} x_m^*\Big\|_{L^{p'}(S;L^2(\O;X^*))}
  \leq C\mu(A_1)^{1/p'}\Big\|\sum_{m=1}^N r_m
  x_m^*\Big\|_{L^2(\O;X^*)}.
\end{equation*}
Since $X$ has non-trivial type, taking the supremum over all
$\sum_{n=1}^N r_n x_n^*\in{\rm Rad}_N(X^*)$ with norm one, it follows that
\begin{equation*}
\begin{split}
   \Big\|\sum_{n=1}^N r_n x_n\Big\|_{L^2(\O;X)}
   &\lesssim\Big\|\sum_{n=1}^N r_n \frac{\one_{A_n}}{\mu(A_1)^{1/p}} x_n\Big\|_{L^p(S;L^2(\O;X))} \\
   &\lesssim\Big\|\sum_{n=1}^N r_n f_n x_n\Big\|_{L^p(S;L^2(\O;X))},
\end{split}
\end{equation*}
where the last estimate used the contraction principle and the fact that \(|f_n|>(2t_0)^{-1/p}\one_{A_n}\gtrsim\mu(A_1)^{-1/p}\one_{A_n}\) by the definition of \(A_n\).

\eqref{it:dualLpContr} : The case $p=\tilde{p}=1$ follows from
\[\begin{aligned}
  LHS\eqref{eq:estLpnormfreverse} &\leq\Big\|\sum_{n=1}^N r_n \int_S|f_n(s)|d\mu(s) x_n \Big\|_{L^2(\O;X)} \\
      &\leq \int_S\Big\|\sum_{n=1}^N r_n |f_n(s)| x_n \Big\|_{L^2(\O;X)}d\mu(s)
        =RHS\eqref{eq:estLpnormfreverse},
\end{aligned}\]
where the first estimate was the contraction principle. For $p>1$, we argue by duality in a similar spirit as in \eqref{it:dualLorContr}: assuming \(\min_{1\leq n\leq N}\|f_n\|_{L^{\tilde{p}}(S)}=1\), choose \(g_n\in L^{\tilde{p}'}(S)\) of at most unit norm so that \(\int_S f_n\cdot g_n\,d\mu=1\) and write
\[
  \sum_{n=1}^N\langle x_n,x_n^*\rangle
  =\E\int_S\Big\langle \sum_{n=1}^N r_n f_n x_n,\sum_{m=1}^N r_m g_m x_m^*\Big\rangle\,d\mu.
\]
Then proceed as in \eqref{it:dualLorContr}, only using
Lemma~\ref{lem:LqrboundedLorentz}\eqref{it:LqContraction} instead of
Lemma~\ref{lem:LqrboundedLorentz}\eqref{it:LorentzContraction}.

\eqref{it:dualLorConverse}: This follows in a similar way as the
corresponding claim in Lemma~\ref{lem:LqrboundedLorentz}.
\end{proof}

\section{Integral operators}

An operator-valued function $T:S\to \calL(X,Y)$ will be called {\em
$X$-strongly measurable} if for all $x\in X$, the $Y$-valued
function $s\mapsto T(s) x$ is strongly measurable.

Let $r\in [1, \infty]$. For an $X$-strongly measurable mapping
$T:S\to \calL(X,Y)$ with $\|T(s) x\|_{L^r(S;Y)} \leq M\|x\|$ and
$f\in L^{r'}(S)$ we will define $T_f\in \calL(X,Y)$ as
\[T_f x = \int_S T(s)x \, f(s) \, d\mu(s).\]
By H\"older's inequality, we have $\|T_{f}\|_{\calL(X,Y)}\leq M
\|f\|_{L^{r'}(S)}$. If $r=1$, then by \cite[Corollary 2.17]{KuWe}
\begin{equation}\label{eq:RboundalwaysL1}
R\big(\{T_f: \|f\|_{L^\infty(S)}\leq 1\}\big)\leq 2 M.
\end{equation}

In the next result we will obtain $R$-boundedness of $\{T_f:
\|f\|_{L^{r'}(S)}\leq 1\}$ for different exponents $r$ under
assumptions on the cotype of $X$ and type of $Y$.

\begin{proposition}\label{prop:intoperator}
Let $X$ and $Y$ be Banach spaces and let $(S,\Sigma, \mu)$ be a
$\sigma$-finite measure space. Let $p_0\in [1, 2]$ and $q_0\in [2,
\infty]$. Assume that $X$ has cotype $q_0$ and $Y$ has type $p_0$.
The following assertions hold:
\begin{enumerate}
\item If $r\in [1, \infty)$ is
such that $\frac1r
>\frac1{p_0}-\frac1{q_0}$. Then there exists a constant
$C=C(r,p_0,q_0,X,Y)$ such that for all $T\in L^r(S;\calL(X,Y))$,
\begin{equation}\label{eq:TfRbounded}
R \big(\{T_f\in \calL(X,Y): \|f\|_{L^{r'}(S)}\leq 1 \}\big) \leq
C\|T\|_{L^r(S;\calL(X,Y))}.
\end{equation}

\item Assume the pair $(p_0, q_0)$ is an element in $\{(1,\infty),(2,\infty),(2,2),(1,2)\}$.
If $r\in [1, \infty]$ is such that $\frac1r =
\frac1{p_0}-\frac1{q_0}$, then there exists a constant $C=C(X,Y)$
such that for all $T\in L^r(S;\calL(X,Y))$ \eqref{eq:TfRbounded}
holds.
\end{enumerate}
\end{proposition}

\begin{remark}\label{rem:nonsep}
Since $\calL(X,Y)$ is usually non-separable, it could happen that
$T:S\to \calL(X,Y)$ is not strongly measurable and therefore not in
$L^r(S;\calL(X,Y))$. However, one can replace the assumption that
$T\in L^r(S;\calL(X,Y))$ by the condition that $T$ is $X$-strongly
measurable and $s\mapsto \|T(s)\|$ is in $L^r(S)$ or is dominated by
a function in $L^r(S)$. This does not affect the assertion in
Proposition \ref{prop:intoperator} and the proof is the same.
\end{remark}

The following will be clear from the proof of Proposition
\ref{prop:intoperator} and Remark \ref{rem:optimal}.

\begin{remark}\label{rem:optimal2}
Let $(A_i,\mathcal{A}_i,\nu_i)$, $i=1, 2$ be a $\sigma$-finite measure space.
Let $1<p_0<p_1\leq 2$ and $2\leq q_1<q_0<\infty$. Let $\widetilde{X}$ be a
Banach space with cotype $q_1$ and $\widetilde{Y}$ be a Banach space with type
$p_1$. If $X = L^{q_0}(A_1;\widetilde{X})$ and $Y =L^{p_0}(A_2;\widetilde{Y})$
, then \eqref{eq:TfRbounded} of Proposition \ref{prop:intoperator} (1) holds
with $\frac1r = \frac{1}{p_0} - \frac{1}{q_0}$.
\end{remark}

\begin{proof}[Proof of Proposition \ref{prop:intoperator}]
(1): \ Let $(f_n)_{n=1}^N$ in $L^{r'}(S)$ be such that
$\sup_n\|f_n\|_{L^{r'}(S)}\leq 1$ and $x_1, \ldots, x_N\in X$.

First assume $p_0>1$ and $q_0<\infty$. Let $p\in (1,p_0)$ and $q\in
(q_0,\infty)$ be such that $\frac1r:=\frac{1}{p} - \frac{1}{q}$, and
hence $\frac{1}{r'}=\frac{1}{p'}+\frac{1}{q}$. Let $g_n =
|f_n|^{r'/q}$ and $h_n = \text{sign}(f_n) |f_n|^{r'/p'}$ for $n=1,
\ldots, N$. Then $\|g_n\|_{L^{q}(S)}\leq 1$,
$\|h_n\|_{L^{p'}(S)}\leq 1$ and $f_n = g_n h_n$ for all $n$. Let
$(y_n^*)_{n=1}^N$ in $Y^*$ be such that $\Big\|\sum_{n=1}^N r_n
y^*_n\Big\|_{L^{p'}(\O;Y^*)}\leq 1$. Then it follows from H\"older's
inequality and $\frac{1}{p} = \frac{1}{r}+\frac{1}{q}$ that
\[\begin{aligned}
\E\Big<\sum_{n=1}^N & r_n T_{f_n} x_n, \sum_{n=1}^N r_n y_n^* \Big>
= \sum_{n=1}^N \lb T_{f_n} x_n, y^*_n\rb
\\ & = \int_S \sum_{n=1}^N \lb g_n(s) T(s) x_n, h_n(s) y^*_n\rb \, d\mu(s)
\\ & = \int_S  \E \Big< T(s)\sum_{n=1}^N r_n g_n(s) x_n, \sum_{n=1}^N r_n h_n(s) y^*_n\Big> \, d\mu(s)
\\ & \leq \Big\| T\sum_{n=1}^N r_n g_n x_n \Big\|_{L^{p}(S\times \O;Y)} \Big\|\sum_{n=1}^N r_n h_n
y^*_n\Big\|_{L^{p'}(S\times \O;Y^*)}
\\ & \leq \|T\|_{L^r(S;\calL(X,Y))} \Big\|\sum_{n=1}^N r_n g_n x_n \Big\|_{L^{q}(S\times \O;X)} \Big\|\sum_{n=1}^N r_n h_n
y^*_n\Big\|_{L^{p'}(S\times \O;Y^*)}.
\end{aligned}\]
Since $X$ has cotype $q_0<q$ it follows from
Lemma~\ref{lem:LqrboundedLorentz}\eqref{it:LqContraction} that
\[\Big\|\sum_{n=1}^N r_n g_n x_n \Big\|_{L^{q}(S\times \O;X)}\leq C_1 \Big\|\sum_{n=1}^N r_n x_n \Big\|_{L^{q}(\O;X)} \]
Since $Y$ has type $p_0$ it follows that $Y^*$ has cotype $p_0'<p'$ (cf.
\cite[Proposition 11.10]{DJT}) and therefore it follows from Lemma~\ref{lem:LqrboundedLorentz}\eqref{it:LqContraction} that
\[\Big\|\sum_{n=1}^N r_n h_n y^*_n\Big\|_{L^{p'}(S\times \O;Y^*)} \leq C_2 \Big\|\sum_{n=1}^N r_n y^*_n\Big\|_{L^{p'}(\O;Y^*)} \leq C_2.\]

We may conclude that
\[\begin{aligned}
\E\Big<\sum_{n=1}^N  r_n T_{f_n} x_n, \sum_{n=1}^N r_n y_n^* \Big>
\leq C_1 C_2 \|T\|_{L^r(S;\calL(X,Y))} \Big\|\sum_{n=1}^N r_n x_n
\Big\|_{L^{q_0}(\O;X)}.
\end{aligned}\]
By assumption $Y$ has non-trivial type, hence ${\rm Rad}_N(Y)^* = {\rm
Rad}_N(Y^*)$ isomorphically (see \eqref{eq:Rad}). Taking the supremum over all
$y_1^*, \ldots, y_N^*\in Y^*$ as above, we obtain that
\[\Big\|\sum_{n=1}^N r_n T_{f_n} x_n\Big\|_{L^{p}(\O;Y)} \lesssim \|T\|_{L^r(S;\calL(X,Y))} \Big\|\sum_{n=1}^N r_n x_n \Big\|_{L^{q}(\O;X)}.\]
The result now follows from \eqref{eq:KKineq}.

If $p_0>1$ and $q_0=\infty$ one can easily adjust the above argument
to obtain the result. In particular, $g_n=1$ for $n=1, \ldots, N$ in
this case.

If $p_0=1$ and $q_0<\infty$, then the duality argument does not hold
since $Y$ only has the trivial type $1$. However, one can argue more
directly in this case. Now $r'>q_0$. By the triangle inequality,
H\"older's inequality and Lemma~\ref{lem:LqrboundedLorentz}, we
obtain that
\[\begin{aligned}
\Big\|\sum_{n=1}^N r_n T_{f_n} x_n\Big\|_{L^{2}(\O;Y)} &\leq \int_S
\Big\|T \sum_{n=1}^N r_n f_n x_n\Big\|_{L^{2}(\O;Y)} \, d\mu(s)
\\ & \leq \|T\|_{L^{r}(S;\calL(X,Y))} \Big\|\sum_{n=1}^N r_n f_n x_n\Big\|_{L^{r'}(S;L^{2}(\O;Y))}
\\ & \leq C \|T\|_{L^{r}(S;\calL(X,Y))} \Big\|\sum_{n=1}^N r_n
x_n\Big\|_{L^{2}(\O;Y)}.
\end{aligned}
\]

(2): \ The case $p_0 = 1$ and $q_0=\infty$ follows from
\eqref{eq:RboundalwaysL1}. The cases ($p_0=2$ and $q_0 = \infty$) and ($p_0
=q_0=2$) follow from Lemma \ref{lem:LqrboundedLorentz} in the same way as
in as (1).

If $p_0=1$ and $q_0 = 2$ then $r=2$ and by the Cauchy-Schwartz inequality and
Lemma \ref{lem:LqrboundedLorentz} we obtain that
\[\begin{aligned}
\Big\|\sum_{n=1}^N r_n T_{f_n} x_n\Big\|_{L^{2}(\O;Y)}& \leq \int_S \Big\|T(s)
\sum_{n=1}^N r_n f_n(s) x_n\Big\|_{L^{2}(\O;X)} \, d\mu(s) \\ & \leq
\|T\|_{L^2(S;\calL(X,Y))} \Big\|\sum_{n=1}^N r_n f_n x_n\Big\|_{L^{2}(S\times
\O;X)}
\\ & \leq C\|T\|_{L^2(S;\calL(X,Y))} \Big\|\sum_{n=1}^N r_n
x_n\Big\|_{L^{2}(S\times \O;X)}.
\end{aligned}\]
\end{proof}

With a certain price to pay, it is possible to relax the norm integrability
condition of Proposition~\ref{prop:intoperator} to the uniform
$L^r$-integrability of the orbits \(s\mapsto T(s)x\). This is reached at the
expense of not being able to exploit the information about the cotype of $X$
but only the type of $Y$, as shown in the following remark. In the example
further below, it is shown that in general the $L^r$-integrability of the
orbits is not sufficient for the full conclusion of
Proposition~\ref{prop:intoperator}.


\begin{remark}\label{rem:intoperatorstrong}
Let $X$ and $Y$ be Banach spaces and let $(S,\Sigma, \mu)$ be a $\sigma$-finite
measure space. Let $p_0\in [1, 2]$. Assume that $Y$ has type $p_0$. The
following assertions hold:
\begin{enumerate}
\item Assume $p_0\in (1, 2)$. If $r\in (1,p_0)$ and $T:S\to \calL(X,Y)$ is such
that
\begin{equation}\label{eq:TstronglyLr}
\|Tx\|_{L^r(S;Y)}\leq C_T \|x\|, \ \ x\in X.
\end{equation}
Then there exists a constant $C=C(r,p_0,Y)$ such that
\begin{equation}\label{eq:TfRboundedstrong}
R \big(\{T_f\in \calL(X,Y): \|f\|_{L^{r'}(S)}\leq 1 \}\big) \leq C C_T.
\end{equation}
\item Assume that $p_0=1$ or $p_0=2$ and that \eqref{eq:TstronglyLr} holds for
$r=p_0$. Then there exists a constant $C=C(Y)$ such that
\eqref{eq:TfRboundedstrong} holds.
\end{enumerate}
\end{remark}

If $Y$ is as in Remark \ref{rem:optimal2}, then \eqref{eq:TfRboundedstrong}
holds for $r=p_0$.

\begin{proof}
(1): \ One can argue as in the proof of Proposition \ref{prop:intoperator} with
$p=r$, $g_n=1$ and $h_n=f_n$. Indeed, we have
\[\begin{aligned}
\E\Big<\sum_{n=1}^N & r_n T_{f_n} x_n, \sum_{n=1}^N r_n y_n^* \Big>
\leq \Big\| T\sum_{n=1}^N r_n x_n \Big\|_{L^{r}(S\times \O;Y)}
\Big\|\sum_{n=1}^N r_n h_n y^*_n\Big\|_{L^{r'}(S\times \O;Y^*)}.
\end{aligned}\]
By the assumption \ref{eq:TstronglyLr} one can estimate
\[\Big\| T\sum_{n=1}^N r_n x_n \Big\|_{L^{p}(S\times \O;Y)}\leq C_T \Big\| \sum_{n=1}^N r_n x_n \Big\|_{L^{p}(S\times \O;X)}.\]
The term $\Big\|\sum_{n=1}^N r_n h_n y^*_n\Big\|_{L^{r'}(S\times \O;Y^*)}$ can
be treated in the same way as in the proof of Proposition
\ref{prop:intoperator}.

(2): \ The case $p_0=1$ follows from \eqref{eq:RboundalwaysL1}. The case
$p_0=2$ can be proved as above.
\end{proof}

In the next example, we will show that even if $X$ has cotype $q$
for some $q\in (2, \infty)$ the result in Remark
\ref{rem:intoperatorstrong} cannot be improved.

\begin{example}
Consider the spaces $X=\ell^q$, $q\in(2,\infty)$, and $Y=\R$, so
that $X$ has cotype $q$ and $Y$ has type $2$. Let $S=\Z_+$ with the
counting measure, and define the $\calL(X,Y)$-valued function $T$ on
$S$ by $T(s)x:=t(s)x(s)$ for some $t:S\to\R$. Then we can make the
following observations:

$T\in L^r(S;\calL(X,Y))$ if and only if $t\in L^r(S)=\ell^r$. The condition~\eqref{eq:TstronglyLr} means $\|tx\|_{r}\lesssim\|x\|_{q}$ (where $tx$ is the pointwise product), which holds if and only if $t\in\ell^u$, $1/u=(1/r-1/q)\vee 0$. Under this condition, the operators $x\mapsto\int_S T(s)f(s)x\,ds$ are well-defined and bounded from $X$ to $Y$ for $f\in L^{r'}(S)=\ell^{r'}$.

Let us then derive a necessary condition for the $R$-boundedness of the mentioned operators. Let $x_i\in X=\ell^q$ be $\alpha(i) e_i$, where $e_i$ is the $i$th standard unit-vector and $\alpha(i)\in\R$. Let $f_i=e_i$. The defining inequality of $R$-boundedness for these functions reduces to
\[
  \|t\alpha\|_2
  \eqsim\E\Big\|\sum_{i=1}^N r_i t(i)\alpha(i)\Big\|
  \lesssim\E\Big\|\sum_{i=1}^N r_i \alpha(i)e_i\|_q
  \eqsim\|\alpha\|_q.
\]
This holds if and only if $t\in\ell^v$, $1/v=1/2-1/q$, which is stronger than~\eqref{eq:TstronglyLr} unless $r\leq 2$. Conversely, the condition~\eqref{eq:TstronglyLr} suffices for $R$-boundedness in this range, as shown in Remark~\ref{rem:intoperatorstrong}.
\end{example}

As a final result in this section we will show how one can use Lemma
\ref{lem:LqrboundedLorentz} to obtain a version of Proposition
\ref{prop:intoperator} with sharp exponents. It may seem artificial
at first sight, but it enables us to obtain a sharp version of
Theorem \ref{thm:main} below. For $f\in L^1(S)+L^\infty(S)$ let
\begin{equation}\label{eq:Lf}
L_{f}(S) = \{g\in L^1(S)+L^\infty(S): \mu(|g|>t) = \mu(|f|>t) \ \text{for all
$t>0$} \}
\end{equation}

\begin{proposition}\label{prop:intoperatorLorentz}
Let $X$ and $Y$ be Banach spaces and let $(S,\Sigma, \mu)$ be a
$\sigma$-finite measure space. Let $p\in [1, 2]$ and $q\in [2,
\infty]$. Assume that $X$ has cotype $q$ and $Y$ has type $p$. If
$r\in [1, \infty]$ is such that $\frac1r =\frac1{p}-\frac1{q}$. Then
there exists a constant $C=C(p,q,X,Y)$ such that for all $f_0\in
L^{r',1}(S)$ and $T\in L^r(S;\calL(X,Y))$,
\begin{equation}\label{eq:TfRboundedLorentz}
R \big(\{T_f\in \calL(X,Y): f\in L_{f_0}(S)\}\big) \leq
C\|T\|_{L^r(S;\calL(X,Y))} \|f_0\|_{L^{r',1}(S)}.
\end{equation}
\end{proposition}

Since each $f\in L_{f_0}(S)$ is also in $L^{r'}(S)$, $T_f\in
\calL(X,Y)$ is well-defined. Note that Remark \ref{rem:nonsep}
applies also here. In the limit cases $p\in \{1,2\}$ and $q\in \{2,
\infty\}$, Proposition \ref{prop:intoperator} (2) yields a stronger
result.

\begin{proof}
Without loss of generality we may assume that
$\|f_0\|_{L^{r,1}(S)}=1$. Let $(f_n)_{n=1}^N$ in $L_{f_0}(S)$ and
$x_1, \ldots, x_N\in X$.

First assume $p>1$ and $q<\infty$. Let $g_n = |f_n|^{r'/q}$ and $h_n
= \text{sign}(f_n) |f_n|^{r'/p'}$ for $n=1, \ldots, N$. Then the
$|g_n|$ have the same distribution as $|f_0|^{r'/q}$ and the $|h_n|$
have the same distribution as $|f_0|^{r'/p'}$, and $f_n = g_n h_n$
for all $n$. Let $(y_n^*)_{n=1}^N$ in $Y^*$ be such that
$\Big\|\sum_{n=1}^N r_n y^*_n\Big\|_{L^{p'}(\O;Y^*)}\leq 1$. Then it
follows from H\"older's inequality and $\frac{1}{p} =
\frac{1}{r}+\frac{1}{q}$ that
\[\begin{aligned}
\E\Big<\sum_{n=1}^N & r_n T_{f_n} x_n, \sum_{n=1}^N r_n y_n^* \Big>
= \sum_{n=1}^N \lb T_{f_n} x_n, y^*_n\rb
\\ & = \int_S \sum_{n=1}^N \lb g_n(s) T(s) x_n, h_n(s) y^*_n\rb \, d\mu(s)
\\ & = \int_S  \E \Big< T(s)\sum_{n=1}^N r_n g_n(s) x_n, \sum_{n=1}^N r_n h_n(s) y^*_n\Big> \, d\mu(s)
\\ & \leq \Big\| T\sum_{n=1}^N r_n g_n x_n \Big\|_{L^{p}(S\times \O;Y)} \Big\|\sum_{n=1}^N r_n h_n
y^*_n\Big\|_{L^{p'}(S\times \O;Y^*)}
\\ & \leq \|T\|_{L^r(S;\calL(X,Y))} \Big\|\sum_{n=1}^N r_n g_n x_n \Big\|_{L^{q}(S\times \O;X)} \Big\|\sum_{n=1}^N r_n h_n
y^*_n\Big\|_{L^{p'}(S\times \O;Y^*)}.
\end{aligned}\]
Since $X$ has cotype $q$ it follows from Lemma
\ref{lem:LqrboundedLorentz} (\ref{it:LorentzContraction}) that
\[\Big\|\sum_{n=1}^N r_n g_n x_n \Big\|_{L^{q}(S\times \O;X)}\leq C_1  \int_0^{\infty}\mu(|f_0|>t^{q/r'})^{1/q} \, dt  \Big\|\sum_{n=1}^N r_n x_n \Big\|_{L^{q}(\O;X)}.\]
It follows from \eqref{eq:Lorentzdistr} and $L^{r',1}\hookrightarrow
L^{r',\frac{r'}{q}}$ that
\[\begin{aligned}
\int_0^{\infty}\mu(|f_0|>t^{q/r'})^{1/q}\, dt &=
\frac{r'}{q}\int_0^{\infty}\mu(|f_0|>t)^{1/q} t^{r'/q} \,\frac{dt}{t} \\ & =
\frac{r'}{q} (r')^{-r'/q} \|f_0\|_{L^{r',\frac{r'}{q}}(S)}^{r'/q}\leq C_{q,r}.
\end{aligned}\]

Since $Y$ has type $p$ it follows that $Y^*$ has cotype $p'$ (cf.
\cite[Proposition 11.10]{DJT}) and therefore it follows from Lemma
\ref{lem:LqrboundedLorentz} (\ref{it:LorentzContraction}) that
\[\Big\|\sum_{n=1}^N r_n h_n y^*_n\Big\|_{L^{p'}(S\times \O;Y^*)} \leq C_2 \int_0^{\infty}\mu(|f_0|>t^{p'/r'})^{1/q}\,dt   \Big\|\sum_{n=1}^N r_n y^*_n\Big\|_{L^{p'}(\O;Y^*)}.\]
As before it holds that
\[\int_0^{\infty}\mu(|f_0|>t^{p'/r'})^{1/q}\,dt\leq C_{p,r} \]
The result now follows with the same duality argument for ${\rm
Rad}_N(Y)$ as in Proposition \ref{prop:intoperator}.

If $p>1$ and $q=\infty$ one can easily adjust the above argument to
obtain the result. In particular, $g_n=1$ for $n=1, \ldots, N$ in
this case. If $p=1$ and $q<\infty$, then one can argue as in
Proposition \ref{prop:intoperator}, but instead of Lemma
\ref{lem:LqrboundedLorentz} (\ref{it:LqContraction}) one has to
apply Lemma \ref{lem:LqrboundedLorentz}
(\ref{it:LorentzContraction}).
If $p=1$ and $q=\infty$ the result follows from Proposition
\ref{prop:intoperator}.
\end{proof}

Similar as in Remark \ref{rem:intoperatorstrong} the following strong version
of Proposition \ref{prop:intoperatorLorentz} holds.

\begin{remark}\label{rem:intoperatorstrongLorentz}
Let $X$ and $Y$ be Banach spaces and let $(S,\Sigma, \mu)$ be a $\sigma$-finite
measure space. Let $p\in (1, 2]$. Assume that $Y$ has type $p$. If $f_0\in
L^{p',1}(S)$ and $T:S\to \calL(X,Y)$ is such that \eqref{eq:TstronglyLr} holds,
then there exists a constant $C=C(p,Y)$ such that
\begin{equation}\label{eq:TfRboundedstrongLorentz}
R \big(\{T_f\in \calL(X,Y): f\in L_{f_0}(S)\}\big) \leq C
C_T\|f_0\|_{L^{p',1}(S)}.
\end{equation}
\end{remark}

\begin{proof}
Similar as in Remark \ref{rem:intoperatorstrong} one can argue as in the proof
of Proposition \ref{prop:intoperatorLorentz} with $p=r$, $g_n=1$ and $h_n=f_n$.

\end{proof}

\section{Besov spaces and $R$-boundedness}


Recall from \cite{We} that for an interval $I=(a,b)$ and $T\in
W^{1,1}(I;\calL(X,Y))$,
\begin{equation}\label{eq:intderivative}
R(T(t)\in \calL(X,Y): t\in (a,b)) \leq \|T(a)\|+
\|T'\|_{L^1(I;\calL(X,Y))}
\end{equation}
In \cite[Theorem 5.1]{GW} this result has been improved under the
assumption that $Y$ has Fourier type $p$. In the next result we
obtain $R$-boundedness for the range of smooth operator-valued
functions under (co)type assumptions on the Banach spaces $X$ and
$Y$. The result below improves \cite[Theorem 5.1]{GW}.
\begin{theorem}\label{thm:main}
Let $X$ and $Y$ be Banach spaces. Let $p\in [1, 2]$ and $q\in [2,
\infty]$. Assume that $X$ has cotype $q$ and $Y$ has type $p$. If
$r\in [1, \infty]$ is such that $\frac1r =\frac1{p}-\frac1{q}$, then
there exists a constant $C=C(p,q,X,Y)$ such that for all $T\in
B^{\frac{d}{r}}_{r,1}(\R^d;\calL(X,Y))$,
\begin{equation}\label{eq:TfRboundedBesov}
R \big(\{T(t)\in \calL(X,Y): t\in \R^d\}\big) \leq C
\|T\|_{B^{\frac{d}{r}}_{r,1}(\R^d;\calL(X,Y))}.
\end{equation}
\end{theorem}

Note that $B^{\frac{d}{r}}_{r,1}(\R^d;\calL(X,Y))\hookrightarrow
BUC(\R^d;\calL(X,Y))$ (the space of bounded, uniformly continuous functions) for all $r\in [1, \infty]$ (cf. \cite[Theorem
2.8.1(c)]{Tr1}).

If $p=q=2$ in (2), then the uniform boundedness of $\{T(t):t\in \R^d\}\subset
\calL(X,Y)$ already implies $R$-boundedness (see \cite[Proposition
1.13]{ArBu}).

%

\begin{proof}
(1): \ Let $Z = \calL(X,Y)$. We may write $T = \sum_{k\geq 0} T_k  =
\sum_{k\geq 0} \sum_{n\geq 0} {\varphi}_n * T_k$, where $T_k =
\varphi_k * T$ and the series converges in $Z$ uniformly in $\R^d$.
Let $\varphi_{-1}=0$. By \cite[Lemma 2.4]{We} we obtain that
\[\begin{aligned}
R\big(T(t)\in Z: t\in \R^d\big) & \leq \sum_{k\geq 0} \sum_{n\geq 0}
R\big(\varphi_n *T_k (t)\in Z: t\in \R^d\big) \\ & = \sum_{k\geq 0}
\sum_{n=k-1}^{k+1} R\big(\varphi_n *T_k (t)\in Z: t\in \R^d\big).
\end{aligned}
\]
Fix $n\in \{0,1,2,\ldots,\}$ and $t\in \R^d$ and define $\varphi_{n,t}\in
{\mathscr S}(\R^d)$ by $\varphi_{n,t}(s) = \varphi_n(t-s)$. Then for all $t\in
\R^d$, $\varphi_{n,t}\in L_{\varphi_n}(\R^d)$, where $L_{\varphi_n}(\R^d)$ is
as in \eqref{eq:Lf} and $\|\varphi_{n}\|_{L^{r',1}(\R^d)}\leq c 2^{\frac{d
n}{r}}$. Indeed, it is elementary to check that $\varphi_n^*(t) = 2^{dn}
\varphi^*(2^{n d}t)$. Therefore,
\[\begin{aligned}
\|\varphi_{n}\|_{L^{r',1}(\R^d)} = 2^{dn} \int_0^\infty t^{\frac{1}{r'}}
\varphi^*(2^{nd}t) \, \frac{dt}{t}  = 2^{\frac{dn}{r}} \int_0^\infty
t^{\frac{1}{r'}} \varphi^*(t) \, \frac{dt}{t} = 2^{\frac{dn}{r}}
\|\varphi\|_{L^{r',1}(\R^d)}
\end{aligned}\]

Letting $T_{k,\varphi_{n,t}}\in Z$ be the integral operator from
Propositions \ref{prop:intoperator} and
\ref{prop:intoperatorLorentz}, it follows that for all $k\geq 0$,
\[\begin{aligned}
R\big(\varphi_{n} *T_k (t)\in Z: t\in \R^d\big) &=
R\big(T_{k,\varphi_{n,t}} \in Z: t\in \R^d\big)   \leq C_1
2^{\frac{d n}{r}} \|T_k\|_{L^r(\R^d;Z)}.
\end{aligned}\]

We may conclude that
\[\begin{aligned}
R\big(T(t)\in Z: t\in \R^d\big) & \leq C_1\sum_{k\geq 0}
\sum_{n=k-1}^{k+1} 2^{\frac{d n}{r}} \|T_k\|_{L^r(\R^d;Z)}
\\ & \leq C_2 \sum_{k\geq 0} 2^{\frac{d k}{r}} \|T_k\|_{L^r(\R^d;Z)}
= C_2 \|T\|_{B^{\frac{d}{r}}_{r,1}(\R^d;Z)}.
\end{aligned}
\]
\end{proof}

%

In \cite[Remark 5.2]{GW} the following result is presented for
operator families which are in a Besov space in the strong sense. If
$Y$ has Fourier type $p$ and $T:\R^d\to \calL(X,Y)$ is such that for
all $x\in X$,
\[
  \|T x\|_{B^{d/p}_{p,1}(\R^d;Y)}\leq C_T \|x\|,
\]
then $\{T(t): t\in \R^d\}$ is $R$-bounded. We will obtain the same conclusion assuming only type $p$.
Notice that many Banach spaces have type $2$, whereas all
Banach spaces with Fourier $2$ are isomorphic to a Hilbert space.

We first need an analogue of Proposition~\ref{prop:intoperator} involving the
space $\gamma(H,Y)$ of $\gamma$-radonifying operators from $H$ to $Y$. See
\cite{NW1} for information on this space. We note that a version of the
following Lemma is also in \cite[Proposition 3.19]{HaKu}, where it is instead
assumed that $Y$ has so-called property $(\alpha)$. Moreover, in \cite[Remark
5.3 and Proof of Proposition 3.19]{HaKu} it is claimed that this assumption can
be relaxed to non-trivial cotype, which is weaker than our assumption below.
However, there seems to be a small confusion there: in \cite[Remark 5.3]{HaKu}
it is observed that non-trivial cotype suffices, thanks to a result
in~\cite{KaiW}, if (a certain Hilbert space) $H_1=\C$, whereas in
\cite[Proposition 3.19]{HaKu} and the following lemma one has the dual
situation: $H_1=H$ is a general Hilbert space and $H_2=\C$. Indeed one could
deduce the following lemma by a standard duality argument from the result
in~\cite{KaiW}, but since a self-contained argument is only slightly longer, we
provide it for completeness:

\begin{lemma}\label{lem:intGamma}
Let $Y$ be a Banach space with non-trivial type and let $H$ be a Hilbert space.
Then there exists a constant
$C$ such that for all $\Psi_n\in\gamma(H,Y)$ and $f_n\in H$,
\[
  \Big\|\sum_{n=1}^N r_n \Psi_n f_n\Big\|_{L^2(\Omega;Y)}
  \leq C\sup_{1\leq n\leq N}\|f_n\|_{H}
   \Big\|\sum_{n=1}^N r_n\Psi_n\Big\|_{L^2(\Omega;\gamma(H,Y))}
\]
\end{lemma}

\begin{proof}
Let $(h_k)_{k=1}^K$ be an orthonormal basis for the span of $(f_n)_{n=1}^N$ in $H$, so that $f_n=\sum_{k=1}^K(h_k,f_n)h_k$. Let further $(r_n)_{n=1}^N$ be a Rademacher sequence on a probability space $\O$, and $(\gamma_k)_{k=1}^K$ a Gaussian sequence on $\O'$. Then
\[\begin{aligned}
  &\Big|\sum_{n=1}^N \langle \Psi_n f_n,y_n^*\rangle\Big|
  =\Big|\sum_{n=1}^N\sum_{k=1}^K \Big<\Psi_n h_k,(h_k,f_n)y_n^*\Big>\Big| \\
  &=\Big|\E\E'\Big<\sum_{k=1}^K \gamma_k\Big(\sum_{n=1}^N r_n \Psi_n\Big)h_k,
         \sum_{m=1}^N r_m\sum_{\ell=1}^K\gamma_{\ell}(h_{\ell},f_m)y_m^*\Big>\Big| \\
  &\leq\Big\|\sum_{k=1}^K \gamma_k\Big(\sum_{n=1}^N r_n \Psi_n\Big)h_k\Big\|_{L^2(\Omega;L^2(\Omega';Y))}
  \\ &\qquad\times\Big\|\sum_{m=1}^N r_m\sum_{\ell=1}^K\gamma_{\ell}(h_{\ell},f_m)y_m^*\Big\|_{L^q(\Omega';L^2(\Omega;Y^*))},
   \qquad q\in[2,\infty).
\end{aligned}\]
The first factor is bounded by
\[
  \Big\|\sum_{n=1}^N r_n \Psi_n\Big\|_{L^2(\Omega;\gamma(H,Y))}
\]
by the definition of the norm in $\gamma(H,Y)$. As for the second, the
non-trivial type of $Y$ implies some non-trivial cotype $q_0\in[2,\infty)$ for
$Y^*$, and then, taking $q\in(q_0,\infty)$ and applying
Lemma~\ref{lem:LqrboundedLorentz}\eqref{it:LqContraction},
\[\begin{aligned}
  &\Big\|\sum_{m=1}^N r_m\sum_{\ell=1}^K\gamma_{\ell}(h_{\ell},f_m)y_m^*\Big\|_{L^q(\Omega';L^2(\Omega;Y^*))} \\
  &\leq C\sup_{1\leq m\leq N}\Big\|\sum_{\ell=1}^K\gamma_{\ell}(h_{\ell},f_m)\Big\|_{L^q(\Omega')}
   \Big\|\sum_{m=1}^N r_m y_m^*\Big\|_{L^2(\Omega;Y^*)}.
\end{aligned}\]
Here $\sum_{\ell=1}^K\gamma_{\ell}(h_{\ell},f_m)$ is a centered Gaussian
variable with variance $\sum_{\ell=1}^K(h_{\ell},f_m)^2=\|f_m\|_{H}^2$, hence
its $L^q$ norm is $c_q\|f_m\|_{H}$ for a constant $c_q$.

The assertion follows by taking the supremum over all $\sum_{m=1}^N r_m y_m^*\in{\rm Rad}_N(Y^*)$ of norm~$1$, using the non-trivial type of~$Y$.
\end{proof}

\begin{proposition}\label{prop:rp} Let $X$ and $Y$ be Banach spaces.
Let $p\in [1, 2]$, and assume that $Y$ has type $p$. If $T:\R^d\to \calL(X,Y)$
satisfies
\begin{equation}\label{eq:TstronglyBesov2}
  \|Tx\|_{B^{d/p}_{p,1}(\R^d;Y)}\leq C_T \|x\|, \ \ x\in X,
\end{equation}
then there exists a constant $C=C(p,Y)$ such that
\begin{equation}\label{eq:TfRboundedBesovstrong2}
  R \big(\{T(t)\in \calL(X,Y): t\in \R^d\}\big) \leq CC_T.
\end{equation}
\end{proposition}



\begin{proof}
If $p=1$, the result follows from \cite[Remark 5.2]{GW}. Let then $p\in (1,2]$.

Fix for the moment $x\in X$ and $k\geq 0$. Let $f_k:\R^d\to Y$ be defined as
$f_k(t) = T_k(t) x = \varphi_k*T(t)x$.
Then by \cite[Theorem 1.1]{KNVW},
\begin{equation}\label{eq:KNVW}
  \|f_k\|_{\g(L^2(\R^d),Y)}
  \leq C\|f_k\|_{B_{p,p}^{d(\frac1p-\frac12)}(\R^d;Y)}
  \leq C2^{d k(\frac{1}{p} - \frac{1}{2})}\|f_k\|_{L^p(\R^d;Y)}.
\end{equation}

Choose $(t_m)_{m=1}^M$ in $\R^d$ and $(x_{m})_{m=1}^M$ in $X$ arbitrarily. Since
$Y$ has type $p>1$ it follows from Lemma~\ref{lem:intGamma} that
\[\begin{aligned}
&\Big\|\sum_{m=1}^M r_m T(t_m)x_m\Big\|_{L^2(\O;Y)}
  \leq\sum_{k\geq 0}\sum_{n=k-1}^{k+1}\Big\|\sum_{m=1}^Mr_m\varphi_n*T_k(t_m)x_m\Big\|_{L^2(\O;Y)} \\
  &=\sum_{k\geq 0}\sum_{n=k-1}^{k+1}\Big\|\sum_{m=1}^M r_m\int_{\R^d}T_k(u)x_m\,\varphi_n(t_m-u)du\Big\|_{L^2(\O;Y)} \\
  &\lesssim\sum_{k\geq 0}\sum_{n=k-1}^{k+1}\sup_{1\leq m\leq M}\|\varphi_n(t_m-\cdot)\|_{L^2(\R^d)}
     \Big\|\sum_{m=1}^M r_m T_k x_m\Big\|_{L^2(\O;\gamma(L^2(\R^d),Y))}
\\  &\eqsim \sum_{k\geq 0} 2^{kd/2}
    \Big\|T_k \Big( \sum_{m=1}^M r_m x_m\Big)\Big\|_{L^2(\O;\gamma(L^2(\R^d),Y))}.
\end{aligned}
\]
Applying \eqref{eq:KNVW} pointwise in $\O$ yields that
\begin{equation}\label{eq:KNVW2}
\Big\|T_k \Big( \sum_{m=1}^M r_m x_m\Big)\Big\|_{\gamma(L^2(\R^d),Y)}\lesssim
2^{kd(\frac1p-\frac12)}\Big\|T_k\Big(\sum_{m=1}^M r_m
x_m\Big)\Big\|_{L^p(\R^d;Y)}.
\end{equation}
Therefore, we obtain from \eqref{eq:KNVW2} and \eqref{eq:KKineq} that
\[\begin{aligned}
\sum_{k\geq 0} & 2^{kd/2}  \Big\|T_k \Big( \sum_{m=1}^M r_m
x_m\Big)\Big\|_{L^2(\O;\gamma(L^2(\R^d),Y))}
\\ & \lesssim \sum_{k\geq 0} 2^{kd/2}  2^{kd(\frac1p-\frac12)} \Big\|T_k\Big(\sum_{m=1}^M r_m x_m\Big)\Big\|_{L^2(\O;L^p(\R^d;Y))}
\\ & \eqsim \int_{\O} \sum_{k\geq 0} 2^{kd/p}\Big\|T_k\Big(\sum_{m=1}^M r_m x_m\Big)\Big\|_{L^p(\R^d;Y)} \, d\P
\\  &=\int_{\Omega}\Big\|T\Big(\sum_{m=1}^M r_m x_m\Big)\Big\|_{B^{d/p}_{p,1}(\R^d;Y)}d\P \\
  &\leq\int_{\Omega} C_T\Big\|\sum_{m=1}^M r_m x_m\Big\|_X d\P
  \leq C_T\Big\|\sum_{m=1}^M r_m x_m\Big\|_{L^2(\Omega;X)}.
\end{aligned}
\]
Putting things together yields the required $R$-boundedness estimate.
\end{proof}



As a consequence of Theorem \ref{thm:main} we have the following two
results. One can similarly derive strong type results from
Proposition~\ref{prop:rp}.

\begin{corollary}\label{cor:integer}
Let $X$ and $Y$ be Banach spaces. Let $p\in [1, 2]$ and $q\in [2,
\infty]$. Assume that $X$ has cotype $q$ and $Y$ has type $p$. Let
$r\in [1, \infty]$ be such that $\frac1r =\frac1{p}-\frac1{q}$. If
there exists an $M$ such that
\begin{equation}\label{eq:TfRboundedSobolev}
\Big(\int_{\R^d}\|D^{\alpha}T\|_{\calL(X,Y)}^r\Big)^{\frac1r}\leq M
\end{equation}
for every $\alpha\in \{0,1,\ldots, d\}^d$ with $|\alpha|\leq
\lfloor\frac{d}{r}\rfloor+1$, then $\{T(t)\in \calL(X,Y): t\in
\R^d\}$ is $R$-bounded.
\end{corollary}

\begin{corollary}\label{cor:holder}
Let $X$ and $Y$ be Banach spaces. Let $p\in [1, 2]$ and $q\in [2, \infty]$.
Assume that $X$ has cotype $q$ and let $Y$ have type $p$. Let $I=(a,b)$ with
$-\infty\leq a<b\leq \infty$. Let $r\in (1, \infty]$ be such that $\frac1r \geq
\frac1{p}-\frac1{q}$. Let $\alpha\in (\frac1r,1)$. If $T\in L^r(\R;\calL(X,Y))$
and there exists an $A$ such that
\begin{equation}\label{eq:Holder}
\|T(s+h)-T(s)\|\leq A|h|^{\alpha} (1+|s|)^{-\alpha}, \ \ s,s+h\in I,
\ h\in I,
\end{equation}
then $\{T(t)\in \calL(X,Y): t\in I\}$ is $R$-bounded by a constant times $A$.
\end{corollary}

Note that in the case that $I$ is bounded, the factor $(1+|s|)^{-\alpha}$ can
be omitted.

\begin{proof}
By taking a worse $p$ or $q$ it suffices to consider the case that
$\frac{1}{r}=\frac1p-\frac1q$. First consider the case that $I=\R$. As in
\cite[Corollary 5.4]{GW} one may check that $T\in \Lambda^{\frac1r}_{r,
1}(\R;\calL(X,Y))$, where the latter is defined in Section \ref{subsec:Besov},
and therefore the result follows from Theorem \ref{thm:main}.

If $I\neq \R$, then one can reduce to the above case by \eqref{eq:gf}.
\end{proof}

\section{Applications}

\subsection{$R$-boundedness of semigroups}
In the next result we will give a sufficient condition for
$R$-boundedness of strongly continuous semigroups restricted to
fractional domain spaces.

\begin{theorem}\label{thm:strongcont}
Let $(T(t))_{t\in \R_+}$ be a strongly continuous semigroup on a Banach space
$X$ with $\|T(t)\|\leq M e^{-\omega t}$ for some $\omega>0$. Assume $X$ has
type $p\in [1, 2]$ and cotype $q\in [2, \infty]$. Let $\alpha>
\frac{1}{r}=\frac1{p}-\frac1{q}$ and let $i_{\alpha}:D((-A)^{\alpha})\to X$ be
the inclusion mapping. Then
\[\{T(t) i_{\alpha}: t\in \R_+\}\subset \calL(D((-A)^{\alpha}), X)\]
is $R$-bounded.
\end{theorem}

\begin{proof}
For $\theta\in (0,1)$ let $X_{\theta} = (X, D(A))_{\theta, \infty}$. Then $x\in
X_{\theta}$ if and only if
\[\|x\|_{X_{\theta}} := \|x\| + \sup_{t\in \R_+} t^{-\theta}\|(T(t) x -x)\|\]
is finite, and this expression defines an equivalent norm on $X_{\theta}$ (cf. \cite[Proposition
3.2.1]{Luninterp}). If we fix $\theta\in (\frac1r,\alpha)$, then we
obtain that
\[\begin{aligned}
\sup_{t\in \R_+} t^{-\alpha} \|T(t) i_{\alpha}
-i_{\alpha}\|_{\calL(D((-A)^{\alpha}),X)} & =
\sup_{\|x\|_{D((-A)^{\alpha})}\leq 1} \sup_{t\in \R_+} t^{-\alpha}
\|T(t) x -x\|_{X} \\ & \leq \sup_{\|x\|_{D((-A)^{\alpha})}\leq 1}
\|x\|_{X_{\theta}}\\ &\lesssim \sup_{\|x\|_{D((-A)^{\alpha})}\leq 1}
\|x\|_{D((-A)^{\alpha})} =1.
\end{aligned}\]
Therefore,
\[\|T(s+h)i_{\alpha}-T(s)i_{\alpha}\|_{\calL(D((-A)^{\alpha}),X)} \lesssim M e^{-\omega s} h^{\alpha}\]
and the result follows from Corollary \ref{cor:holder}.
\end{proof}

The result in Theorem \ref{thm:strongcont} is quite sharp as follows
from the next example. An application of Theorem
\ref{thm:strongcont} will be given in Theorem \ref{thm:applSCP}.

For $\alpha\in \R$ and $p\in [1, \infty]$, let $H^{\alpha,p}(\R)$ be the
Bessel-potential spaces (cf. \cite[2.3.3]{Tr1}).
\begin{example}
Let $p\in [1, \infty)$. Let $(T(t))_{t\in \R}$ be the
left-translation group on $X=L^p(\R)$ with generator $A =
\frac{d}{dx}$. Then for all $\alpha\in (|\frac{1}{p}-\frac12|,1)$
and $M\in \R_+$,
\begin{equation}\label{eq:translation}
\{T(t) i_{\alpha} : t\in [-M,M]\}\subset \calL(H^{\alpha,p}(\R),
L^p(\R)),
\end{equation}
is $R$-bounded, where $i_{\alpha}:H^{\alpha,p}(\R)\to L^p(\R)$
denotes the embedding.

On the other hand, for $\alpha\in (0, |\frac{1}{p}-\frac12|)$ and
$M=1$, the family \eqref{eq:translation} is not $R$-bounded.
\end{example}

\begin{proof}
Note that $L^p(\R)$ and $H^{\alpha,p}$ have type $p\wedge 2$ and cotype $p\vee
2$. Therefore, for $\alpha>|\frac{1}{p}-\frac12|$ the $R$-boundedness of
\[\{e^{-t} T(t) i_{\alpha} : t\in \R_+\}\subset
\calL(H^{\alpha,p}(\R), L^p(\R))\] follows from Theorem \ref{thm:strongcont}.
Therefore, we obtain from the Kahane-contraction principle that
\[\{T(t) i_{\alpha} : t\in [0,M]\}\subset \calL(H^{\alpha,p}(\R),
L^p(\R))\]
is $R$-bounded. Since a similar argument works for $T(-t)$, the
$R$-boundedness of \eqref{eq:translation} follows from the fact that the union
of two $R$-bounded sets is again $R$-bounded.

For the converse, let $\psi\in C^{\infty}(\R)\setminus \{0\}$ be
such that $\text{supp}(\psi)\subset (0,1)$. For $c\in (0,\infty)$
let $\psi_c(t) = \psi(ct)$. Then $(-A)^{\alpha} \psi_c =
c^{\alpha}[(-A)^{\alpha}\psi]_c$.
Fix an integer $N$ and let $f_n = f_0:=
\psi_N$ for all $n$. Then $f_0$ has support in
$(0,1/N)$ and $\|f_0\|_{L^p(\R)}^p=N^{-1}\|\psi\|_{L^p(\R)}^p$.

There holds, on the one hand,
\[\begin{aligned}
  \Big\|\sum_{n=1}^N &r_n T(n/N) i_{\alpha} f_n \Big\|_{L^2(\O;X)}^p
  = \Big\|\sum_{n=1}^N r_n f_0(\cdot+n/N) \Big\|_{L^2(\O;X)}^p \\
  & = \sum_{n=1}^N \|f_0(\cdot+ n/N)\|_{L^p(\R)}^p
  = N\|f_0\|_{L^p(\R)}^p = \|\psi\|_{L^p(\R)}^p,
\end{aligned}\]
and on the other hand,
\[\begin{aligned}
  \Big\|\sum_{n=1}^N & r_n f_n \Big\|_{L^2(\O;D((-A)^{\alpha}))}^p
  = N^{\frac{p}{2}} \|f_0\|_{D((-A)^{\alpha})}^p
  = N^{\frac{p}{2}} \big(\|f_0\|_{L^p(\R)}^p+\|(-A)^{\alpha}f_0\|_{L^p(\R)}^p\big) \\
  &= N^{\frac{p}{2}} \big(\|\psi_N\|_{L^p(\R)}^p+\|N^{\alpha}[(-A)^{\alpha}\psi]_N\|_{L^p(\R)}^p\big)  \\
  &= N^{-1+\frac{p}{2}}\big(\|\psi\|_{L^p(\R)}^p+N^{\alpha p}\|(-A)^{\alpha}\psi]\|_{L^p(\R)}^p\big).
\end{aligned}\]
Therefore, if $\tau:=\{T(t) i_{\alpha}:t\in [-1,1]\}$ is
$R$-bounded, then it follows that there exists a constant $C$ such
that
\[
  1 \leq C N^{-\frac1p+\frac{1}{2}+\alpha}.
\]
Letting $N$ tend to infinity, this implies that $\alpha\geq\frac1p-\frac12$, i.e., $\tau$ can only be $R$-bounded in this range.

We still have to prove that the $R$-boundedness also implies
$\alpha\geq\frac12-\frac1p$. This can be proved by duality. If
$\{T(t)\in\calL(H^{\alpha,p}(\R),L^p(\R)): t\in [-1,1]\}$ is $R$-bounded, then
$\{T^*(t)\in \calL(L^{p'}(\R),H^{-\alpha,p'}(\R)): t\in [-1,1]\}$ is
$R$-bounded as well. It follows that
$\{(1-A)^{-\alpha}T^*(t)\in \calL(L^{p'}(\R),L^{p'}(\R)):t\in [-1,1]\}$ is
$R$-bounded. This implies that $\{T^*(t)\in
\calL(H^{\alpha,p'}(\R),L^{p'}(\R)):t\in [-1,1]\}$ is $R$-bounded. According to
the first part of the proof this implies that $\alpha\geq \frac{1}{p'} -
\frac12 = \frac12-\frac{1}{p}$.
\end{proof}

\subsection{Stochastic Cauchy problems}

We apply Theorem \ref{thm:strongcont} to stochastic equations with additive
Brownian noise. We refer the reader to \cite{NW1} for details on stochastic
Cauchy problems, stochastic integration and $\g$-radonifying operators. Let
$(\O,\F,\P)$ be a probability space. Let $H$ be a separable Hilbert space and
let $W_H$ be a cylindrical Wiener process. Recall from \cite{NW1} that for an
operator-valued function $\Phi:[0,t]\to \calL(H,E)$ which belongs to
$\g(L^2(0,t;H),X)$ (the space of $\gamma$-radonifying operators from $L^2(0,t;H)$ to $X$) we have
\[\Big\|\int_0^t \Phi(s) \, d W_H(s)\Big\|_{L^2(\O;X)} = \|\Phi\|_{\g(L^2(0,t;H),X)}.\]

On a real Banach space $X$ we consider the following equation.
\begin{equation}\tag{SE}\label{SE}
\left\{\begin{aligned}
dU(t) & = AU(t)\,dt + B(t)dW_H(t), \qquad t\in \R_+,\\
 U(0) & = x,
\end{aligned}
\right.
\end{equation}
Here $A$ is the generator of a strongly continuous semigroup $(T(t))_{t\in
\R_+}$, $B:\R_+\to \calL(H,E)$ and $x\in X$. We say that a strongly measurable
process $U:\R_+\times\O\to X$ is a mild solution of \eqref{SE} if for all $t\in
\R_+$, almost surely we have
\[U(t) = T(t) x + \int_0^t T(t-s) B(s) \, d W_H(s).\]

In general \eqref{SE} does not have a solution (cf. \cite[Example 7.3]{NW1}).
In the case when $B(t) = B\in \g(H,X)$ is constant, there are some sufficient
conditions for existence. Indeed, if $X$ has type $2$ or $(T(t))_{t\in \R_+}$
is an analytic semigroup, then \eqref{SE} always has a unique mild solution and
it has a version with continuous paths (see \cite[Corollary 3.4]{VeZi} and
\cite{DevNWe} respectively). In the next result we prove such an existence and
regularity result under assumptions on the noise in terms of the type and
cotype of $X$.

\begin{theorem}\label{thm:applSCP}
Assume $X$ has type $p\in [1, 2]$ and cotype $q\in [2,\infty]$. Let $w\in \R$
be such that $\lim_{t\to \infty} e^{w t}T(t)=0$. Let $\alpha>\frac{1}{p} -
\frac{1}{q}$ and $B\in \g(L^2(\R_+;H),D((w-A)^{\alpha}))$. Then \eqref{SE} has
a unique mild solution $U$. Moreover, if there exists an $\e>0$ such that for
all $M\in \R_+$,
\begin{equation}\label{eq:epsgamma}
\sup_{t\in [0,M]} \|s\mapsto (t-s)^{-\e}
B(s)\|_{\g(L^2(0,t;H),D((w-A)^{\alpha})}<\infty
\end{equation}
then $U$ has a version with continuous paths.
\end{theorem}

In particular we note that if $B(t) = B\in \g(H,D((w-A)^{\alpha}))$ is constant
then for all $\e\in (0,\frac12)$
\[\|s\mapsto (t-s)^{-\e} B(s)\|_{\g(L^2(0,t;H),D((w-A)^{\alpha}))} = (1-2\e)^{-1} t^{\frac12-\e} \|B\|_{\g(H,D((w-A)^{\alpha}))}.\]

\begin{remark}
Here is a sufficient condition for \eqref{eq:epsgamma}: there is an $s\in
(2,\infty)$ such that for all $M\in \R_+$,
\[B\in B^{\frac{1}{p}-\frac12}_{s,p}(0,M;D((w-A)^{\alpha})).\]
Indeed, it follows from \cite[Lemma 3.3]{NVW3} that \eqref{eq:epsgamma} holds
for all $\e\in (0,\frac12-\frac1s)$.
\end{remark}

\begin{proof}
Assume that \eqref{eq:epsgamma} holds for some $\e\in [0,\frac12)$. In the case
$\e=0$ we will show existence of a solution, and in the other case we show that
the solution has a version with continuous paths.

By Theorem \ref{thm:strongcont}, $\{e^{w t} T(t) i_{\alpha}\in
\calL(D((w-A)^{\alpha}),X):t\geq 0\}$ is $R$-bounded. It follows that for fixed
$M>0$,
\[\{T(t) i_{\alpha}\in \calL(D((w-A)^{\alpha}),X):t\in [0,M]\}\]
is $R$-bounded by some constant $C$. Therefore, by \cite{KaWe}
(see also \cite[Theorem 9.14]{vNisem}), the function $s\mapsto T(s)i_{\alpha}$ acts as a
multiplier between the spaces $\g(L^2(0,t;H),D((w-A)^{\alpha}))$ and
$\g(L^2(0,t;H),X)$, and we conclude that
\[\begin{aligned}
  \sup_{t\in [0,M]}&\|s\mapsto (t-s)^{-\e}T(t-s) B(s)\|_{\g(L^2(0,t;H),X)} \\
  & \leq C\sup_{t\in [0,M]}\|s\mapsto(t-s)^{-\e}B(s)\|_{\g(L^2(0,t;H),D((w-A)^\alpha))}
  <\infty.
\end{aligned}\]
Now the result follows from \cite[Proposition 3.1 and Theorem 3.3]{VeZi}.
\end{proof}

\subsection{$R$-boundedness of evolution families\label{subsec:AT}}

In the next application we obtain $R$-boundedness of an evolution
family generated by a family $(A(t))_{t\in [0,T]}$ of unbounded
operators which satisfy the conditions (AT) of Acquistapace and
Terreni (see \cite{AT2}). For $\phi\in(0,\pi]$, we define the sector
\[  \Sigma(\phi):=\{0\}\cup\{\lambda\in\C\setminus\{0\}:
    |\arg(\lambda)|<\phi\}. \]

The condition (AT) is said to be satisfied if the following two
requirements hold:
\begin{enumerate}
\item[(AT1)]  \label{AT1}
The $A(t)$ are linear operators on a Banach space $ E $ and there
are constants $K \ge 0$, and $ \phi \in (\frac{\pi}{2},\pi)$ such
that $\Sigma(\phi) \subset \varrho(A(t))$ and for all $ \lambda \in
\Sigma(\phi)$ and $t\in [0,T]$,
\[ \| R(\lambda, A(t)) \| \le \frac{K}{1+ |\lambda|}.\]
\item[(AT2)]  \label{AT2}
There are constants $L \ge 0$ and $\mu, \nu \in (0,1]$ with $\mu +
\nu >1$ such that for all $\lambda \in \Sigma(\phi,0)$ and $s,t \in
[0,T]$,
\[ \| A(t)R(\lambda,
A(t))(A(t)^{-1}- A(s)^{-1})\| \le L |t-s|^\mu (|\lambda|+1)^{-\nu}.
\]
\end{enumerate}
Under these assumptions there exists a unique strongly continuous
evolution family $(P(t,s))_{0\leq s\leq t\leq T}$ in $\calL(X)$ such
that $\frac{\partial P(t,s)}{\partial t} = A(t) P(t,s)$ for $0\leq
s<t\leq T$. Moreover, $\|A(t)P(t,s)\|\leq C(t-s)^{-1}$.

For analytic semigroup generators one has that for all $\e>0$ and
$T\in [0,\infty)$, $\{t^{\e}S(t)\in \calL(X): t\in [0,T]\}$ is
$R$-bounded. This easily follows from \eqref{eq:intderivative}. This
may be generalized to evolution families $(P(t,s))_{0\leq s\leq
t\leq T}$, where $(A(t))_{t\in [0,T]}$ satisfies the (AT)
conditions. Indeed, then by the same reasoning we obtain that for
all $\alpha>0$,
\[\sup_{s\in [0,T]} R\big(\{(t-s)^{\alpha}P(t,s)\in \calL(X): t\in [s,T]\}\big)<\infty.\]
This argument does not hold if one considers the $R$-bound with
respect to $s\in [0,t]$ instead of $t\in [s,T]$. This is due to the
fact that
\begin{equation}\label{eq:diffPs}
\Big\|\frac{\partial P(t,s)}{\partial s}\Big\|\leq
C(t-s)^{-1}
\end{equation}
might not be true. The $R$-boundedness with respect to $s\in [0,t]$
has applications for instance in the study of non-autonomous
stochastic Cauchy problems (see \cite{VeZi}). We also note that
\eqref{eq:diffPs} does hold if $(A(t)^*)_{t\in [0,T]}$ satisfies the
(AT)-conditions (see~\cite{AT3}).

Recall from \cite[Theorem 2.3]{Ya} that for all $\theta\in (0,\mu)$,
\begin{equation}\label{eq:Yatheta}
\|P(t,s)(-A(s))^{\theta}\| \leq C(t-s)^{-\theta}, \ 0\leq s<t\leq T.
\end{equation}
Due to this inequality one might expect that under assumptions on
$\mu$, one can still obtain a fractional version of
\eqref{eq:diffPs}. This is indeed the case and in the next theorem
we will give conditions under which the $R$-boundedness with respect
to $s\in [0,t]$ holds.

The authors are grateful to Roland Schnaubelt for showing them the following
result.

\begin{proposition}\label{prop:analytic}
Assume (AT). Then for all $\theta\in (0,\mu)$ there exists a constant $C$ such
that for all $0\leq s\leq t\leq T$,
\begin{equation}\label{eq:analyticityest}
\|(-A(t))^{-\theta}(P(t,s)-I)\|\leq C (t-s)^{\theta}.
\end{equation}
\end{proposition}

\begin{proof}
First let $\theta\in (1-\nu,\mu)$. By \cite[equation (A.5)]{MSchn} we can write
\[(-A(t))^{-\theta}(P(t,s)-I)  =  g(t,s) + \int_s^{t} (-A(t))^{-\theta}P(t,\tau) (-A(\tau))^{\theta} h(\tau,s) d\tau,\]
where
\[g(t,s) = (-A(t))^{-\theta} (e^{(t-s)A(s)}-I),\]
\[h(t,s) = (-A(t))^{1-\theta} \big[(-A(s))^{-1} - (-A(t))^{-1} \big] A(s) e^{(t-s)A(s)}.\]
We may write
\[\begin{aligned}
  \|g(t,s)\| \leq &\|((-A(t))^{-\theta} - (-A(s))^{-\theta}) (e^{(t-s)A(s)}-I)\| \\
    &+ \| (-A(s))^{-\theta} (e^{(t-s)A(s)}-I)\|.
\end{aligned}\]
By \cite[equation (2.10)]{Schn}
\[\|((-A(t))^{-\theta} - (-A(s))^{-\theta}) (e^{(t-s)A(s)}-I)\|\lesssim |t-s|^{\mu}\lesssim (t-s)^{\theta}.\]
For the other term it is clear that
\[\| (-A(s))^{-\theta} (e^{(t-s)A(s)}-I)\| \leq \int_0^{t-s} \|(-A(s))^{1-\theta} e^{\tau A(s)} \| \, d\tau \lesssim (t-s)^{\theta}.\]
This shows that $\|g(t,s)\|\lesssim (t-s)^{\theta}$. By \cite[equation
(2.2)]{Ya90} we obtain that
\[\|h(t,s)\| \lesssim (t-s)^{\mu-1}.\]
Since by \cite[Lemma A.1]{MSchn} $V(t,s) = (-A(t))^{-\theta}P(t,\tau)
(-A(\tau))^{\theta}$ is uniformly bounded, it follows that
\[\int_s^{t} \|V(t,\tau) h(\tau,s)\| d\tau\lesssim (t-s)^{\mu}\lesssim (t-s)^{\theta}.\]
We may conclude \eqref{eq:analyticityest} for the special choice of $\theta$.

For general $\theta\in (0,\mu)$ choose $\e>0$ so small that $\mu-\e>1-\nu$.
Then by interpolation with $\theta/(\mu-\e)$ it follows that
\[\begin{aligned}
\|& (-A(t))^{-\theta}(P(t,s)-I)\| \\ & \lesssim
\|(-A(t))^{-(\mu-\e)}(P(t,s)-I)\|^{\theta/(\mu-\e)}
\|P(t,s)-I\|^{1-\theta/(\mu-\e)} \lesssim (t-s)^{\theta}.
\end{aligned}\]
\end{proof}

\begin{theorem}
Let $X$ be a Banach space with type $p\in [1,2]$ and cotype $q\in [2, \infty]$.
Assume (AT) with
\[\mu>\frac{1}{p} - \frac{1}{q}.\]
Then for all $\e>0$,
\[\sup_{t\in [0,T]} R\big(\{(t-s)^{\e}P(t,s)\in \calL(X): s\in
[0,t]\}\big)<\infty .\]
\end{theorem}

As a consequence one obtains a version of \cite[Corollary 4.5]{VeZi}
without assuming $\|\frac{\partial P(t,s)}{\partial s}\|\leq
C(t-s)^{-1}$.

\begin{proof}
Choose $\theta\in (\frac{1}{p} - \frac{1}{q},\mu)$. Let $r\in (1,
\infty)$ be such that $\theta>\frac1r \geq \frac{1}{p} -
\frac{1}{q}$ and $\e>\theta-\frac1r$. Fix $t\in [0,T]$. We will
apply Theorem \ref{thm:main}, with the equivalent norm explained in
Section \ref{subsec:Besov}, to the function $f:[0,t]\to \calL(X)$
defined by $f(s) = (t-s)^{\e}P(t,s)$.

Let $h\in (0,t)$. By the triangle inequality we can write
\[\begin{aligned}
\|&(t-s-h)^{\e} P(t,s+h) - (t-s)^{\e}P(t,s)\| \\ & \leq
\|(t-s-h)^{\e}(P(t,s+h) -P(t,s))\| + |(t-s-h)^{\e}- (t-s)^{\e}|
\|P(t,s)\|.
\end{aligned}\]

Since the main point is dealing with small $\e$, we may assume that
$(\e-1)r<-1$. Then
\[\begin{aligned}
  &\int_0^{t-h}|(t-s)^{\e}-(t-s-h)^{\e}|^r\, ds
   =\Big(\int_h^{2h}+\int_{2h}^t\Big)|u^{\e}-(u-h)^{\e}|^r\,du \\
  &\leq\int_h^{2h}u^{\e r}\,du+\int_{2h}^t |h\e(u-h)^{\e-1}|^r\,du
   \leq h(2h)^{\e r}+(h\e)^r\int_h^{\infty}u^{(\e-1)r}\,du \\
  &\lesssim_{\e,r} h^{1+\e r}+h^{r}h^{1+(\e-1)r}\eqsim h^{1+\e r}.
\end{aligned}\]
For the other part it follows from \eqref{eq:Yatheta} and
Proposition \ref{prop:analytic} that for all $\theta<\mu$
\[\begin{aligned}
\|P(t,s+h) -P(t,s)\| &= \|P(t,s+h)(-A(s+h))^{\theta}
(-A(s+h))^{-\theta} (I- P(s+h,s))\|
\\ &\leq C(t-s-h)^{-\theta}\|(-A(s+h))^{-\theta} (P(s+h,s)-I)\|
\\ &\leq C(t-s-h)^{-\theta} h^{\theta}.
\end{aligned}\]

We conclude that
\[\begin{aligned}
  \Big(\int_0^{t-h} \|&(t-s-h)^{\e} P(t,s+h) - (t-s)^{\e}P(t,s)\|^r \, ds\Big)^{\frac1r} \\
  & \lesssim h^{\e+\frac1r} + h^{\theta} \Big(
     \int_0^{t-h} (t-s-h)^{(\e-\theta)r} \, ds \Big)^{\frac1r} \\
  & \eqsim h^{\e+\frac1r} + h^{\theta} \lesssim h^{\theta}
\end{aligned}\]
where we used $\e>\theta-\frac1r$. Similar results hold for $h<0$.

It follows that $\varrho_r(f,\tau)\lesssim \tau^{\theta}$, for
$\tau\in (0,1)$, where $\varrho_r$ is defined as in Section
\ref{subsec:Besov}. Since $\theta>\frac1r$ we can conclude that
$f\in \Lambda^{\frac1r}_{r, 1}(0,t;\mathcal{B}(X))$. Now the result
follows from \eqref{eq:gf}, the norm equivalence of
$\Lambda^{\frac1r}_{r, 1}(\R;\mathcal{B}(X))$ and $B^{\frac1r}_{r,
1}(\R;\mathcal{B}(X))$, and Theorem \ref{thm:main}.
\end{proof}

{\em Acknowledgment} -- The authors thank S. Kwapie\'n for the
helpful comments on Lemma~\ref{lem:LqrboundedLorentz} which
eventually led us to parts \eqref{it:LorentzContraction} and
\eqref{it:LorentzConverse} of that Lemma. The authors thank R.
Schnaubelt for showing them Proposition~\ref{prop:analytic}.

\end{document}